\documentclass[a4paper,fleqn]{cas-sc}

\usepackage[authoryear,longnamesfirst]{natbib}
\usepackage{amsmath,amssymb,amsthm,mathtools}
\usepackage{bm}
\usepackage{graphicx}
\graphicspath{{Graphics/}}
\usepackage{booktabs}
\usepackage{array}
\usepackage{multirow}
\usepackage{subcaption}
\usepackage{placeins}

\theoremstyle{plain}
\newtheorem{theorem}{Theorem}
\newtheorem{lemma}[theorem]{Lemma}
\newtheorem{assumption}[theorem]{Assumption}
\newtheorem{corollary}[theorem]{Corollary}
\theoremstyle{remark}
\newtheorem{remark}[theorem]{Remark}

\newcommand{\R}{\mathbb{R}}
\newcommand{\Sphere}{\mathbb{S}^2}

\newcommand{\norm}[1]{\left\lVert #1 \right\rVert}
\newcommand{\ip}[2]{\left(#1,#2\right)}
\newcommand{\dd}{\,\mathrm{d}}
\newcommand{\eps}{\varepsilon}
\newcommand{\dt}{\Delta t}
\newcommand{\gradS}{\nabla_{\!\Sphere}}
\newcommand{\LapS}{\Delta_{\!\Sphere}}

\newcommand{\CapSpinodal}{Spinodal decomposition on $\Sphere$ at $T=2$
($\eps=0.01$, $r=4$, $p=5$, semi-implicit, $\dt = 5\times10^{-4}$).
Two views of the order parameter $u_h$.  Sharp interfaces of width
$O(\eps)$ separate phase-separated domains.  The colour scale is
symmetric and clipped at the extreme values of $u_h$.}

\newcommand{\CapHpSnapshots}{Spinodal decomposition at $T=2$
($\eps=0.01$, $\dt = 5\times10^{-4}$, semi-implicit) across the three
mesh configurations of Table~\ref{tab:meshes}. From left to right:
r4$p$3, r3$p$6, and r2$p$12 (23\,042 DOFs each).
A common symmetric colour scale is used for all three panels.}

\newcommand{\CapStructure}{Structure diagnostics for the semi-implicit
scheme ($r=4$, $p=5$, $\dt=5\times10^{-4}$, $\eps=0.01$, $S=600$),
recorded at every step of the Wolfram Language~15.0.1 rerun.
Left: the absolute quadrature mass $|m_Q^n|$ of the zero-mass initial-value
problem remains at $O(10^{-13})$. Right: the full history of the truncated
discrete energy $\mathcal E_h(u^n)$ of \eqref{eq:Eh}; the inset resolves
$0\le t\le0.1$, including the rise from $165.4$ at $t=0.060$ to $222.8$
at $t=0.075$. The modified Lyapunov functional $E_h^n$ of \eqref{eq:Lyap}
decreases throughout this rerun and equals $124.7$ when $\mathcal E_h=85.3$
at $t=0.3$. Table~\ref{tab:energy_mono} gives the original-run increase
counts, with the rerun counts stated in the accompanying text.}

\newcommand{\CapTorusPlain}{Non-viscous Cahn--Hilliard coarsening on the
torus at $t=0$, $0.05$, $0.20$, $0.50$, and $1.50$ (left to right).
The truncated double-well potential of \eqref{eq:trunc} is used, with $\eps=0.025$,
$\dt=5\times10^{-4}$, $S=12$, $p=5$, and $T=1.5$. All panels share one symmetric
colour scale.}

\newcommand{\CapTorusVoigt}{Viscous R{\"a}tz--Voigt coarsening on the
torus at $t=0$, $0.1$, $0.25$, $0.4$, $0.55$, $0.7$, $0.85$, and $1.0$
(left to right), computed with the stabilized SAV-BDF2 integrator of
Section~\ref{sec:voigt-test}.  Parameters are $\gamma=0.05$, $\alpha=\nu=1$,
$S=360$, $\dt=10^{-3}$, and $p=3$.  The Voigt double-well places the two phases
at $u=0$ and $u=1$. All panels share one colour scale.}

\newcommand{\CapGenusTwoGrid}{Cahn--Hilliard coarsening on the genus-2
implicit surface \eqref{eq:phi-genus2} at
$t=0$, $0.05$, $0.1$, $0.15$, $0.2$, $0.3$, $0.5$, $0.7$, and $1$
(left to right, top to bottom).  Parameters: $p=3$, $h=0.3$,
$\eps=0.1$, $\dt=2\times10^{-3}$, semi-implicit BDF2.  The rendering
evaluates the order-$p$ field on a $6\times$ subdivided, exactly
reprojected curved geometry.}

\newcommand{\CapGenusTwoHero}{Final state ($t=1$) of the genus-2 run of
Figure~\ref{fig:genus2_grid}. The Ginzburg--Landau energy is
$\mathcal E_h = 11.877$.  The two phases have organized into a single
interface pair wrapping both handles.  Left: tilted view from above.
Right: tilted view from below, both at $\sim\!60^\circ$ elevation so the
two handle holes are visible.}

\newcommand{\CapGenusThreeGrid}{Cahn--Hilliard coarsening on the genus-3
implicit surface \eqref{eq:phi-genus3} at
$t=0$, $0.05$, $0.1$, $0.15$, $0.2$, $0.3$, $0.5$, $0.7$, and $1$
(left to right, top to bottom).  Parameters: $p=3$, $h=0.3$,
$\eps=0.1$, $\dt=2\times10^{-3}$, semi-implicit BDF2. Rendering as in
Figure~\ref{fig:genus2_grid}.}

\newcommand{\CapGenusThreeHero}{Final state ($t=1$) of the genus-3 run of
Figure~\ref{fig:genus3_grid}. The Ginzburg--Landau energy is
$\mathcal E_h = 24.677$.  The banded morphology threads all three
handles. Coarsening is still ongoing at $T=1$.  Left: tilted view from
above. Right: tilted view from below, both at $\sim\!60^\circ$ elevation
so the three handle holes are visible.}

\begin{document}
\let\WriteBookmarks\relax
\def\floatpagepagefraction{1}
\def\textpagefraction{.001}

\shorttitle{Cahn--Hilliard on Closed Surfaces}
\shortauthors{H. Hakula, A. Rasila, X. Xu}

\title[mode=title]{Semi-Implicit and Convex Splitting Schemes for the Surface Cahn--Hilliard Equation on Curved Closed Surfaces}

\author[1]{Harri Hakula}
\ead{harri.hakula@aalto.fi}
\affiliation[1]{organization={Department of Mathematics and Systems Analysis, Aalto University},
  city={Espoo},
  country={Finland}}

\author[2,4]{Antti Rasila}
\cormark[1]
\cortext[cor1]{Corresponding author}
\ead{antti.rasila@iki.fi}
\ead{antti.rasila@gtiit.edu.cn}
\affiliation[2]{organization={Department of Mathematics with Computer Science, Guangdong Technion -- Israel Institute of Technology},
  city={Shantou},
  country={China}}

\author[3,4]{Xinpeng Xu}
\ead{xu.xinpeng@gtiit.edu.cn}
\affiliation[3]{organization={Department of Physics and Guangdong Key Laboratory of Materials and Technologies for Energy Conversion, Guangdong Technion -- Israel Institute of Technology},
  city={Shantou},
  country={China}}

  \affiliation[4]{organization={Technion -- Israel Institute of Technology},
  city={Haifa},
  country={Israel}}

\begin{abstract}
The surface Cahn--Hilliard equation describes phase separation and coarsening of a conserved composition field on fixed curved surfaces. We study a mixed $hp$-surface finite element discretization with curvilinear high-order elements, primarily on the sphere, with the stabilized semi-implicit BDF2 integrator and a stabilized convex-splitting variant solved with the same factorization. For the semi-implicit scheme with a quadratic-growth potential and a truncated nonlinear bulk chemical-potential term, we prove exact discrete mass conservation, energy decay for the start-up steps, an exact per-step energy balance for the BDF2 phase, and a conditional monotone law for a modified Lyapunov functional.
A manufactured solution with analytic forcing verifies second order in time and the expected $h$- and $p$-convergence. For $\eps = 0.01$, the sphere runs conserve mass to linear-solve accuracy, the modified Lyapunov functional is monotone at every step at the two larger time steps, and the balance closes to roundoff. Production time steps exceed the fastest linear time scale: the stabilized scheme then flattens the dispersion relation, and a post-onset refinement study finds order-one time-step errors at the two largest steps and a few percent, relative to its reference, at a step of $2\times10^{-5}$. The large-step runs provide qualitative phase-separated configurations and tests of discrete conservation and modified-energy behavior. Their trajectories and final configurations depend on the time step. At fixed mesh and time step, the semi-implicit implementation is 4--6 times faster than the convex-splitting variant, with different potentials outside $[-1,1]$.
Torus tests include controlled equal-area bands, whose diffuse energies reproduce sharp-interface perimeters within $1\%$, deterministic coarsening, and a viscous R{\"a}tz--Voigt model. Genus-two and genus-three calculations on implicit surfaces conserve mass to roundoff, and their unstable-mode counts agree with area-based Weyl estimates within $1\%$.
\end{abstract}

\begin{highlights}
\item The mixed $hp$-FEM discretization preserves mass exactly at the fully discrete level.
\item Sphere growth rates and torus band energies are checked against analytical values.
\item An exact per-step energy balance closes to roundoff along stabilized BDF2 runs.
\item Equal-dimension $hp$ tests isolate assembly, solve, and conditioning trade-offs.
\end{highlights}


\maketitle

\section{Introduction}

Time-dependent partial differential equations on closed surfaces $\Gamma \subset \mathbb{R}^3$ arise when transport, phase separation, or orientational order is confined to an interface. Applications include multicomponent lipid membranes, vesicles, and membrane-bound active layers~\citep{BaumgartHessWebb2003,JulicherLipowsky1993,CallanJonesSorreBassereau2011,ParthasarathyYuGroves2006}. The surface Cahn--Hilliard equation is a conserved gradient-flow model for lateral demixing and coarsening on such interfaces~\citep{CahnHilliard1958,RatzVoigt2006,MarenduzzoOrlandini2013}.

The passive model studied here couples composition to geometry through the surface metric in $\nabla_\Gamma$, $\Delta_\Gamma$, and the surface measure. These local metric effects can change growth rates and interface lengths and bias domain organization. Topology supplies a separate global constraint by determining the homology classes available to phase boundaries. Explicit composition--curvature coupling can introduce curvature-dependent chemical preferences on prescribed surfaces~\citep{MarenduzzoOrlandini2013,JulicherLipowsky1993,CallanJonesSorreBassereau2011}.

A sphere and an ellipsoid both have genus zero but different metrics,
whereas a torus has genus one and two inequivalent noncontractible
cycles. These distinctions occur in closed-membrane models and in
surface-PDE benchmarks with higher-genus or nonuniform
geometries~\citep{MichaletBensimon1995,BouzarMenasMueller2015}. Surface
finite element and related geometry-flexible methods address such
settings~\citep{DziukElliott2013,RatzVoigt2006,DuJuTian2011,
ElliottRanner2015,ZhangWangRong2011,Ratz2016Benchmark}. Building on
$hp$-discretizations for Laplace--Beltrami problems on closed
surfaces~\citep{HakulaRasila2025, hrz}, we study discrete mass conservation,
modified-energy diagnostics for semi-implicit BDF2 stepping, and
$p$-refinement for thin-interface resolution on curved surfaces.

The Cahn--Hilliard equation models phase separation through the $H^{-1}$-gradient flow of the Ginzburg--Landau energy ~\citep{CahnHilliard1958}. On closed surfaces, it provides a standard diffuse-interface model for phase separation constrained to a manifold, see, for example, ~\citep{DuJuTian2011,RatzVoigt2006,ElliottRanner2015,LiKimWang2017}. On the unit sphere $\Sphere$, the equation reads
\begin{align}
\partial_t u &= \LapS \mu,
\label{eq:ch1}
\\
\mu &= -\eps \LapS u + \frac{1}{\eps} f(u),
\label{eq:ch2}
\end{align}
where $u$ is the order parameter, $\mu$ is the chemical potential, and $f(u)=u^3-u$ is the derivative of the standard double-well potential. The contribution $\eps^{-1}f(u)$ is the nonlinear bulk chemical-potential term. We use ``reaction term'' and ``reaction vector'' as numerical terminology for this local nonlinearity and its assembled finite element vector. The mobility has been nondimensionalised to unity throughout, except in the viscous R{\"a}tz--Voigt benchmark of Section~\ref{sec:voigt-test}, where it is stated explicitly. Here
\[
F(u)=\frac14(u^2-1)^2.
\]

The stabilized semi-implicit BDF2 scheme with vector reaction extrapolation $2f(u^n)-f(u^{n-1})$ and second-difference stabilization goes back to the stabilized schemes of \citet{XuTang2006} and \citet{ShenYang2010}, appears in the form used here in \citep[Eq.~(6.22a)]{Shen2011}, and has been applied to the Cahn--Hilliard equation on the sphere by \citet{ZhangWangRong2011} with a prolate spectral-element discretization. For surface Cahn--Hilliard problems, \citet{LiKimWang2017} give a linear, second-order, mass-conserving, unconditionally energy-stable Crank--Nicolson scheme with a stabilized splitting on triangulated surfaces, \citet{DuJuTian2011} analyze finite element approximations, and \citet{ElliottRanner2015} analyze evolving-surface finite elements.

%
\subsection{Main results}\label{sec:main-results}

The mixed $hp$ formulation on curved surface elements provides the
following analytical and numerical results.
The identity $Kc=0$ yields exact mass conservation for both schemes and
a mass invariant for the convex-splitting iterates up to linear-solve roundoff
(Section~\ref{sec:mass}). For the truncated semi-implicit scheme, we prove
start-up energy decay, derive an exact BDF2 energy balance, and obtain
decay of a modified Lyapunov functional under a bound on the extrapolation
defect (Section~\ref{sec:energy}). Trajectory diagnostics verify the balance
to roundoff and test the defect bound directly (Section~\ref{sec:defect}).

A manufactured solution demonstrates second-order temporal convergence
and the expected $h$- and exponential $p$-convergence. A resolved
single-mode test recovers the analytical growth rate. At production
parameters, stabilization flattens the dispersion relation, and post-onset
refinement reveals order-one errors at the two largest time steps relative
to a finer-step reference (Section~\ref{sec:production-accuracy}).

The convex-splitting solver reuses the semi-implicit factorization for
its nonlinear corrections. In the sphere production runs, semi-implicit
stepping is $4.2$--$5.8$ times faster at fixed mesh and time step
(Table~\ref{tab:dt_comparison}). This compares implementation costs for
schemes whose potentials differ outside $[-1,1]$.
At equal dimension and resolution scale, the $hp$ comparison separates
assembly and solve costs and identifies the intermediate r3$p$6
configuration as having the lowest stepping cost (Table~\ref{tab:hp_comparison}).

Controlled torus bands reproduce analytical interface perimeters within
$1\%$. The deterministic torus run changes from hole-dominant to
tube-dominant Fourier content and ends with a connected phase domain
classified as wrapping both periodic directions. On implicit
genus-two and genus-three surfaces, mass is conserved to roundoff and
unstable-mode counts agree with area-based Weyl estimates within $1\%$
(Section~\ref{sec:geometry-tests}).

\subsection{Organization}

Section~\ref{sec:SI} defines the mixed spatial discretization and the
stabilized semi-implicit scheme. Section~\ref{sec:cs} introduces the
convex-splitting variant and its fixed-matrix solver.
Section~\ref{sec:structurelaws} recalls the continuous mass and energy
identities. Section~\ref{sec:mass} proves exact discrete mass conservation,
and Section~\ref{sec:energy} derives start-up energy decay, the exact BDF2
energy balance, and the conditional Lyapunov law.
Section~\ref{sec:experiments} presents accuracy verification, tests of the
discrete structure laws, time-step and cost studies, and geometry and
topology benchmarks. Section~\ref{sec:conclusion} summarizes the findings.

\section{Semi-Implicit Discretization}\label{sec:SI}

The spatial discretization is a mixed $hp$ surface finite element method. Time is advanced by $n_0$ stabilized backward-Euler start-up steps followed by a stabilized semi-implicit BDF2 method. The convergence tests and the torus and implicit-surface runs use $n_0=1$. The $\eps=0.01$ sphere production runs use $n_0=2$, an implementation choice whose effect is recorded in Section~\ref{sec:production-setup}. This section defines the matrices, truncated reaction term, quadrature-defined potential, and linear systems used in the analysis and experiments.

Let $V_h \subset H^1(\Sphere)$ be the $hp$-finite element space. Denote by $M,K\in\R^{N_h\times N_h}$ the mass and stiffness matrices. For coefficient vectors $v,w\in\R^{N_h}$, write
\[
\ip{v}{w}_M := v^T M w,
\qquad
\ip{v}{w}_K := v^T K w,
\qquad
\norm{v}_M^2 := v^T M v,
\qquad
\norm{v}_K^2 := v^T K v.
\]
We assume that constants belong to $V_h$ and denote by $c\in\R^{N_h}$ the coefficient vector of the discrete constant $1_h\equiv 1$, so that
\begin{equation}
Kc = 0.
\label{eq:Kc0}
\end{equation}
The identity \eqref{eq:Kc0} is the algebraic form of the closed-surface fact that the Laplace--Beltrami operator has constants in its kernel. Exact mass conservation below follows from this identity and the first block row of the mixed scheme.

For the semi-implicit scheme, the reaction term is replaced by the quadratic-growth truncation
\begin{equation}
\tilde f(s)=
\begin{cases}
2(s+1), & s<-1,\\
s^3-s, & |s|\le 1,\\
2(s-1), & s>1,
\end{cases}
\qquad
\tilde F(s)=
\begin{cases}
(s+1)^2, & s<-1,\\
\tfrac14(s^2-1)^2, & |s|\le 1,\\
(s-1)^2, & s>1.
\end{cases}
\label{eq:trunc}
\end{equation}
Then $\tilde F'=\tilde f$ and $\tilde f$ is globally Lipschitz with
\begin{equation}
L := \operatorname*{ess\,sup}_{s\in\R} |\tilde f'(s)| = 2.
\label{eq:L2}
\end{equation}
The equality follows from $\tilde f'(s)=2$ for $|s|>1$ and $\tilde f'(s)=3s^2-1$ for $|s|\le1$, whose maximum absolute value on $[-1,1]$ is attained at $s=\pm1$. This truncation follows \citet{CondetteMelcherSuli2011}. The semi-implicit scheme and the structure-preserving statements below are formulated for this truncated nonlinearity. When the nodal and quadrature-point values remain in $[-1,1]$, the scheme coincides with the standard quartic Cahn--Hilliard discretization. Overshoots activate the truncated discrete model. To match the discrete model and the numerical diagnostics, define the discrete potential with the same quadrature rule used to assemble the reaction vector:
\begin{equation}
\mathcal R_h(u) := \sum_{q\in Q_h} w_q\, \tilde F\bigl(u_h(x_q)\bigr).
\label{eq:Rh}
\end{equation}
On the quartic branch $|u_h|\le1$, the degree-$4p+1$ production quadrature used below integrates the degree-$4p$ polynomial $\tilde F(u_h)$ exactly in the reference coordinates of each element, up to the curved-element surface Jacobian. When the truncation is active, $\mathcal R_h$ is the quadrature-defined discrete potential used by the scheme. Since $\tilde F$ is piecewise polynomial, this quantity may differ from the continuous integral of $\tilde F(u_h)$ across cells intersecting the kink at $u_h=\pm1$.
Then
\begin{equation}
D\mathcal R_h(u)[v] = v^T R_{\tilde f}(u)
\qquad \text{for all } u,v\in\R^{N_h},
\label{eq:Rh-grad}
\end{equation}
and the truncated discrete energy monitored in the experiments is
\begin{equation}
\mathcal E_h(u) := \frac{\eps}{2}\norm{u}_K^2 + \frac1{\eps}\mathcal R_h(u).
\label{eq:Eh}
\end{equation}

The time discretization uses $n_0\in\{1,2\}$ stabilized semi-implicit backward-Euler start-up steps, for $k=1,\dots,n_0$,
\begin{align}
\frac{1}{\dt} M\bigl(u^k-u^{k-1}\bigr) + K\mu^k &= 0,
\label{eq:ie-1}
\\
M\mu^k &= \eps K u^k + \frac{S}{\eps} M\bigl(u^k-u^{k-1}\bigr) + \frac{1}{\eps} R_{\tilde f}(u^{k-1}).
\label{eq:ie-2}
\end{align}
For every $n\ge n_0$, the subsequent BDF2 step is
\begin{align}
\frac{1}{\dt} M\Bigl(\frac32 u^{n+1} - 2u^n + \frac12 u^{n-1}\Bigr) + K\mu^{n+1} &= 0,
\label{eq:bdf2-1}
\\
M\mu^{n+1} &= \eps K u^{n+1}
+ \frac{S}{\eps} M\bigl(u^{n+1} - u^{\mathrm{ext},n}\bigr)
+ \frac{1}{\eps} R_{\tilde f}^{\mathrm{ext},n},
\label{eq:bdf2-2}
\end{align}
with
\begin{equation}
u^{\mathrm{ext},n} := 2u^n-u^{n-1},
\qquad
R_{\tilde f}^{\mathrm{ext},n} := 2R_{\tilde f}(u^n)-R_{\tilde f}(u^{n-1}).
\label{eq:ext}
\end{equation}
The semi-implicit linear systems are uniquely solvable. Eliminating $\mu$ gives Schur complements
\[
\frac1{\dt}M + K M^{-1}\Bigl(\eps K+\frac{S}{\eps}M\Bigr)
\quad\text{and}\quad
\frac{3}{2\dt}M + K M^{-1}\Bigl(\eps K+\frac{S}{\eps}M\Bigr),
\]
which are symmetric positive definite because $M$ is positive definite, $K$ is positive semidefinite, and $\dt,\eps,S>0$.

\section{Convex-Splitting Variant}\label{sec:cs}
For comparison we also consider a stabilized convex-splitting BDF2 variant on the same spatial discretization \citep{Eyre1998}. This section specifies the splitting, the nonlinear equations, the fixed-mass affine spaces on which they are posed, and the fixed-matrix correction used to solve them. The description matches the solver that produced the convex-splitting rows of Table~\ref{tab:dt_comparison}.

The variant splits the untruncated quartic double well,
\begin{equation}
F(s) = F_c(s) + F_e(s),
\qquad
F_c(s) := \tfrac14 s^4 + \tfrac14,
\qquad
F_e(s) := -\tfrac12 s^2,
\label{eq:cs-split}
\end{equation}
so that $f_c(s):=F_c'(s)=s^3$ is monotone and $f_e(s):=F_e'(s)=-s$ is linear. The cubic reaction vector $R_c(u)$ is assembled from $f_c$ with the reaction rule of Section~\ref{sec:spatial}. Since $f_e$ is linear, the expansive reaction vector is taken as $R_e(u)=-Mu$, for which the vector and state extrapolations coincide. The implicit treatment of the full cubic relies on its monotonicity. The convex-splitting range diagnostics are available at the terminal states, which lie within $[-1.07,1.15]$ (Table~\ref{tab:dt_comparison}). The semi-implicit trajectories at the same steps reach $|u_h|\approx3$ transiently (Section~\ref{sec:structure}). For $|u_h|>1$, the cubic $f_c$ grows faster than the truncated $\tilde f$ but remains monotone, so the two variants solve different discrete models there. The truncated diagnostics $\mathcal R_h$ and $\mathcal E_h$ of \eqref{eq:Rh}--\eqref{eq:Eh} serve as common monitors on the convex-splitting states.

The scheme keeps the stabilization term of \eqref{eq:bdf2-2} and treats the untruncated cubic implicitly. With $u^{\mathrm{ext},n}=2u^n-u^{n-1}$ as in \eqref{eq:ext}, the first start-up step is
\begin{align}
\frac{1}{\dt} M\bigl(u^1-u^0\bigr) + K\mu^1 &= 0,
\label{eq:cs-start-1}
\\
M\mu^1 &= \eps K u^1 + \frac{S}{\eps}M\bigl(u^1-u^0\bigr)
+ \frac{1}{\eps}\bigl(R_c(u^1) - Mu^0\bigr),
\label{eq:cs-start-2}
\end{align}
the second start-up step, as implemented, keeps the backward-Euler mass term but extrapolates the explicit part from $(u^0,u^1)$,
\begin{align}
\frac{1}{\dt} M\bigl(u^2-u^1\bigr) + K\mu^2 &= 0,
\label{eq:cs-start2-1}
\\
M\mu^2 &= \eps K u^2 + \frac{S}{\eps}M\bigl(u^2-u^{\mathrm{ext},1}\bigr)
+ \frac{1}{\eps}\bigl(R_c(u^2) - Mu^{\mathrm{ext},1}\bigr),
\label{eq:cs-start2-2}
\end{align}
and for $n\ge2$ the BDF2 step is
\begin{align}
\frac{1}{\dt} M\Bigl(\tfrac32 u^{n+1} - 2u^n + \tfrac12 u^{n-1}\Bigr) + K\mu^{n+1} &= 0,
\label{eq:cs-1}
\\
M\mu^{n+1} &= \eps K u^{n+1} + \frac{S}{\eps}M\bigl(u^{n+1}-u^{\mathrm{ext},n}\bigr)
+ \frac{1}{\eps}\bigl(R_c(u^{n+1}) - Mu^{\mathrm{ext},n}\bigr).
\label{eq:cs-2}
\end{align}
Equivalently, the implicit part is $f_c(s)+Ss$ and the extrapolated part is $(S+1)s$, which is the form in which the solver assembles the equations. Both start-up steps are first-order consistent. Equations \eqref{eq:cs-start-1}--\eqref{eq:cs-2} are nonlinear in the new state through $R_c$.

The first block rows coincide with those of the semi-implicit scheme. The start-up steps are posed on the affine space $\{v:\ip{c}{v}_M=\ip{c}{u^0}_M\}$, where the first row is range-compatible with $K$. For the BDF2 step, if the two previous states have common mass $m^0$, the first row is range-compatible with $K$ on the affine space
\[
\{v:\ip{c}{v}_M=m^0\},
\]
since this condition is equivalent to
\[
c^T M\Bigl(\tfrac32u^{n+1}-2u^n+\tfrac12u^{n-1}\Bigr)=0,
\]
under the induction hypothesis $\ip{c}{u^{n-1}}_M=\ip{c}{u^n}_M=m^0$. On these prescribed-mass affine subspaces the nonlinear contribution $R_c(u)+SMu$ is strictly monotone because $F_c$ is convex and $M$ is positive definite. Existence follows from the corresponding finite-dimensional strictly convex minimization problems, and uniqueness from the strict monotonicity supplied by the backward-Euler or BDF2 mass term together with $\eps K$ and the monotone reaction contribution.

\paragraph{Fixed-matrix corrections.}
The nonlinear equations are solved by corrections with the fixed block matrix
\begin{equation}
P_{\mathrm{SI}}
:=
\begin{bmatrix}
\tfrac{3}{2\dt} M & K \\[2pt]
\eps K + \dfrac{S}{\eps} M & -M
\end{bmatrix},
\label{eq:PSI}
\end{equation}
which is the semi-implicit BDF2 block matrix, factored once and reused (the backward-Euler block matrix, with $\tfrac{3}{2\dt}$ replaced by $\dt^{-1}$, serves the start-up steps). The residuals are
$F_1=\dt^{-1}M(\tfrac32u^{n+1}-2u^n+\tfrac12u^{n-1})+K\mu^{n+1}$ and
$F_2=\eps K u^{n+1}+\eps^{-1}\bigl(R_c(u^{n+1})+SMu^{n+1}-(S+1)Mu^{\mathrm{ext},n}\bigr)-M\mu^{n+1}$,
the initial value is $u^{n+1,0}=u^{\mathrm{ext},n}$ with $\mu^{n+1,0}$ obtained from the second row at $u^{n+1,0}$, and each correction solves $P_{\mathrm{SI}}(\Delta u^k,\Delta\mu^k)=-(F_1^k,F_2^k)$. The first row of the correction,
\[
\frac{3}{2\dt}M\Delta u^k+K\Delta\mu^k=-F_1^k,
\]
is the property used for the iteration-level mass invariant in Theorem~\ref{thm:cs-mass}.

The Jacobian of $(F_1,F_2)$ differs from $P_{\mathrm{SI}}$ only in the $(2,1)$ block, by $\eps^{-1}J_c(u)$, where $J_c(u)=\partial R_c/\partial u$ is the symmetric positive semidefinite matrix assembled from the weight $3u_h^2$ with the common rule. The corrections are therefore a quasi-Newton iteration whose linearized error propagation on a Laplace--Beltrami eigenmode with $K\phi=\lambda M\phi$ and a locally constant weight $3\bar u^2$ has the factor
\begin{equation}
\rho(\lambda)=\frac{(3\bar u^2/\eps)\,\lambda}{\tfrac{3}{2\dt}+\eps\lambda^2+(S/\eps)\lambda}
\le \frac{3\bar u^2}{S}.
\label{eq:cs-contraction}
\end{equation}
With $S=600$ the bound is $6.6\times10^{-3}$ at $|u_h|=1.15$ and $4.5\times10^{-2}$ at $|u_h|=3$. At these contraction rates, reducing the scaled residual from $10^{3}$--$10^{4}$ to $10^{-10}$ would require six to ten corrections. The observed $3$--$6$ corrections per step are consistent in scale with this estimate (Section~\ref{sec:solver_comparison}). The first correction is more effective because the extrapolation error is concentrated in low modes, where $\rho(\lambda)$ is small. Equation~\eqref{eq:cs-contraction} describes modal contraction under a constant-weight linearization. Omitting the implicit $SMu^{n+1}$ term makes $P_{\mathrm{SI}}$ differ from the Jacobian by $(S/\eps-3\bar u^2/\eps)M$ in the $(2,1)$ block. The factor then becomes $\approx(S-3\bar u^2)/S$ on every resolved mode, and the corrections stall (Section~\ref{sec:solver_comparison}).

\begin{remark}
Relative to the classical first-order Eyre scheme, \eqref{eq:cs-1}--\eqref{eq:cs-2} uses BDF2 time stepping, extrapolates the expansive part, and carries the same second-difference stabilization as the semi-implicit scheme. The implicit cubic term supplies the monotone nonlinear contribution. The schemes share the second-difference stabilization but differ in nonlinear time treatment and, outside $[-1,1]$, in the free-energy potential being discretized. Section~\ref{sec:mass} records the corresponding mass invariant, Section~\ref{sec:solver_comparison} the numerical diagnostics, and Section~\ref{sec:production-accuracy} the accuracy comparison.
\end{remark}

\section{Continuous Structure Laws}\label{sec:structurelaws}
Before turning to the discrete conservation and dissipation statements, we recall the two continuous identities that the numerical schemes are meant to respect. Closedness of the surface gives zero boundary flux, so the Cahn--Hilliard flow conserves the total mass and dissipates the free energy through the chemical potential. These elementary identities provide the template for the discrete results in the next two sections.

For smooth solutions of \eqref{eq:ch1}--\eqref{eq:ch2}, the continuous energy is
\begin{equation}
\mathcal E(u) = \int_{\Sphere} \frac{\eps}{2} |\gradS u|^2 + \frac{1}{\eps} F(u)\,\dd S.
\label{eq:cont-energy}
\end{equation}

\begin{lemma}[Continuous mass and energy laws]
Let $(u,\mu)$ be a smooth solution of \eqref{eq:ch1}--\eqref{eq:ch2} on the closed surface $\Sphere$. Then
\begin{align}
\frac{\dd}{\dd t}\int_{\Sphere} u\,\dd S &= 0,
\label{eq:cont-mass}
\\
\frac{\dd}{\dd t}\mathcal E(u) &= - \int_{\Sphere} |\gradS \mu|^2\,\dd S \le 0.
\label{eq:cont-diss}
\end{align}
\end{lemma}

\begin{proof}
Integrating \eqref{eq:ch1} over $\Sphere$ and using closedness of $\Sphere$ gives \eqref{eq:cont-mass}. For \eqref{eq:cont-diss}, differentiate \eqref{eq:cont-energy} in time, substitute \eqref{eq:ch2}, and integrate by parts:
\[
\frac{\dd}{\dd t}\mathcal E(u)
= \int_{\Sphere} \mu\,\partial_t u\,\dd S
= \int_{\Sphere} \mu\,\LapS \mu\,\dd S
= - \int_{\Sphere} |\gradS \mu|^2\,\dd S.
\qedhere
\]
\end{proof}

\section{Exact Discrete Mass Conservation}\label{sec:mass}
We now verify the discrete analogue of the continuous mass law. The argument uses the discrete constant, the identity $Kc=0$, and the first block row. It applies to the semi-implicit scheme, to the convex-splitting nonlinear solve, and to its fixed-preconditioned iterates.

\begin{theorem}\label{thm:mass}
The semi-implicit discretization \eqref{eq:ie-1}--\eqref{eq:bdf2-2} preserves mass exactly:
\[
\ip{c}{u^n}_M = \ip{c}{u^0}_M
\qquad \text{for all } n\ge 0.
\]
\end{theorem}

\begin{proof}
Left-multiply \eqref{eq:ie-1} by $c^T$. Using \eqref{eq:Kc0} yields
\[
\ip{c}{u^k-u^{k-1}}_M = 0,\qquad k=1,\dots,n_0,
\]
hence $\ip{c}{u^k}_M=\ip{c}{u^0}_M$ for every start-up step.

For the BDF2 step, left-multiply \eqref{eq:bdf2-1} by $c^T$ and use \eqref{eq:Kc0} again:
\[
\frac32 \ip{c}{u^{n+1}}_M - 2 \ip{c}{u^n}_M + \frac12 \ip{c}{u^{n-1}}_M = 0,
\qquad n\ge n_0.
\]
Writing $m^n := \ip{c}{u^n}_M$, the recurrence
\[
\frac32 m^{n+1} - 2m^n + \frac12 m^{n-1} = 0
\]
with $m^{n_0}=m^{n_0-1}=m^0$ implies $m^{n_0+1}=m^{n_0}$, and then $m^{n+1}=m^n$ for all $n\ge n_0$ by induction.
\end{proof}

\begin{theorem}[Mass invariant for convex-splitting solves]\label{thm:cs-mass}
Assume that the convex-splitting nonlinear systems are uniquely solvable on the prescribed-mass affine subspaces described in Section~\ref{sec:cs}. Then the start-up steps \eqref{eq:cs-start-1}--\eqref{eq:cs-start2-2} and the BDF2 step \eqref{eq:cs-1}--\eqref{eq:cs-2} preserve the discrete mass:
\[
\ip{c}{u^n}_M=\ip{c}{u^0}_M
\qquad \text{for all } n\ge0.
\]
Moreover, if the start-up nonlinear systems are solved by the fixed-matrix corrections of Section~\ref{sec:cs} initialized with $u^{1,0}=u^0$ and $u^{2,0}=u^{\mathrm{ext},1}$, and if the BDF2 nonlinear system is initialized with $u^{n+1,0}=u^{\mathrm{ext},n}$, then every nonlinear iterate has the initial mass $m^0$, up to linear-solve roundoff.
\end{theorem}

\begin{proof}
For the exact nonlinear start-up and BDF2 solutions on the prescribed-mass affine subspaces, mass preservation follows by left-multiplying the first block row, \eqref{eq:cs-start-1}, \eqref{eq:cs-start2-1}, or \eqref{eq:cs-1}, by $c^T$ and using $Kc=0$, exactly as in Theorem~\ref{thm:mass}.

For the start-up iterations, suppose the current iterate $u^{k,j}$ has mass $m^0$. The first residual block
\[
F_1^{k,j}=\dt^{-1}M(u^{k,j}-u^{k-1})+K\mu^{k,j}
\]
satisfies $c^T F_1^{k,j}=0$, and the first block row of the corresponding correction gives $\ip{c}{\Delta u^j}_M=0$. Since the initializations $u^{1,0}=u^0$ and $u^{2,0}=2u^1-u^0$ have mass $m^0$, all start-up iterates preserve mass.

For the BDF2 iteration, assume by induction that $\ip{c}{u^{n-1}}_M=\ip{c}{u^n}_M=m^0$. Then
\[
\ip{c}{u^{\mathrm{ext},n}}_M
=2\ip{c}{u^n}_M-\ip{c}{u^{n-1}}_M
=m^0,
\]
so the initial iterate has the correct mass. Suppose the current iterate $u^{n+1,k}$ has mass $m^0$. The first residual block
\[
F_1^k=\dt^{-1}M\Bigl(\tfrac32u^{n+1,k}-2u^n+\tfrac12u^{n-1}\Bigr)+K\mu^{n+1,k}
\]
satisfies $c^T F_1^k=0$. The correction solves
\[
\frac{3}{2\dt}M\Delta u^k+K\Delta\mu^k=-F_1^k
\]
in the first block row. Left-multiplication by $c^T$ gives $\ip{c}{\Delta u^k}_M=0$, hence the next iterate also has mass $m^0$. The converged nonlinear iterate therefore inherits the same mass invariant. The mass argument depends only on the first block row.
\end{proof}

\section{Discrete Energy Laws}\label{sec:energy}
The two-step nature of BDF2 enters the structure analysis through the energy estimate. Under the start-up condition of Lemma~\ref{lem:startup}, each backward-Euler start-up step dissipates the truncated discrete energy $\mathcal E_h$. Under Assumption~\ref{ass:ext}, the BDF2 phase dissipates a modified Lyapunov functional $E_h^n$ containing the BDF2 $G$-norm contribution and the stabilization increment.

The proof is organized as a sequence of short algebraic estimates. The $G$-identity handles the BDF2 gradient term, the quadrature-based Lipschitz and Taylor bounds handle the truncated reaction, and the start-up lemma treats the backward-Euler steps separately. Corollary~\ref{cor:balance} derives the exact balance \eqref{eq:balance} directly from the BDF2 update. The extrapolation defect in that balance can have either sign. Assumption~\ref{ass:ext} bounds the defect, giving the conditional decay law of Theorem~\ref{thm:lyap}. The numerical diagnostics test the exact balance directly.

\begin{lemma}[BDF2 $G$-identity]\label{lem:Gid}
For every symmetric positive semidefinite matrix $A$ and every $a,b,d\in\R^m$,
\begin{equation}
\ip{\frac32 a - 2b + \frac12 d}{a}_A
= \frac14\Bigl(
\norm{a}_A^2 + \norm{2a-b}_A^2
- \norm{b}_A^2 - \norm{2b-d}_A^2
+ \norm{a-2b+d}_A^2
\Bigr).
\label{eq:Gid}
\end{equation}
\end{lemma}

\begin{proof}
Let
\[
G=\frac14
\begin{bmatrix}
5 & -2\\
-2 & 1
\end{bmatrix},
\qquad
\norm{(v,w)}_{G,A}^2
:=
\begin{bmatrix}v\\ w\end{bmatrix}^{\!T}
\begin{bmatrix}
G_{11}A & G_{12}A\\
G_{21}A & G_{22}A
\end{bmatrix}
\begin{bmatrix}v\\ w\end{bmatrix}.
\]
Then
\[
\norm{(v,w)}_{G,A}^2
=\frac14\bigl(\norm{v}_A^2+\norm{2v-w}_A^2\bigr),
\]
and hence
\[
\norm{(a,b)}_{G,A}^2-\norm{(b,d)}_{G,A}^2
=\frac14\Bigl(\norm{a}_A^2+\norm{2a-b}_A^2-\norm{b}_A^2-\norm{2b-d}_A^2\Bigr).
\]
Adding $\frac14\norm{a-2b+d}_A^2$ and expanding the last term gives
$\frac32\norm{a}_A^2-2\ip{a}{b}_A+\frac12\ip{a}{d}_A$, which is the left-hand side of \eqref{eq:Gid}.
\end{proof}

\begin{lemma}[Discrete Lipschitz bound]\label{lem:Lip}
Assume that the mass matrix $M$ and the reaction vector $R_{\tilde f}$ are assembled with a common positive quadrature rule $\{(x_q,w_q)\}_{q\in Q_h}$, so that the quadrature representation
\begin{equation}
z^T M w = \sum_{q\in Q_h} w_q\, z_h(x_q)\, w_h(x_q),
\qquad
z^T R_{\tilde f}(u) = \sum_{q\in Q_h} w_q\, z_h(x_q)\, \tilde f\bigl(u_h(x_q)\bigr)
\label{eq:quad-rep}
\end{equation}
holds for all $z,w,u\in\R^{N_h}$ with associated finite element functions $z_h,w_h,u_h\in V_h$. Then, with $L$ as in \eqref{eq:L2},
\begin{equation}
\norm{M^{-1}\bigl(R_{\tilde f}(u)-R_{\tilde f}(v)\bigr)}_M
\le L\, \norm{u-v}_M
\qquad \text{for all } u,v\in\R^{N_h}.
\label{eq:LipM}
\end{equation}
\end{lemma}

\begin{proof}
For any $z\in\R^{N_h}$, \eqref{eq:quad-rep} gives
\[
z^T\bigl(R_{\tilde f}(u)-R_{\tilde f}(v)\bigr)
= \sum_{q\in Q_h} w_q\, z_h(x_q)\bigl(\tilde f(u_h(x_q))-\tilde f(v_h(x_q))\bigr).
\]
Since $\tilde f$ is globally Lipschitz with constant $L$ (cf.\ \eqref{eq:L2}), the Cauchy--Schwarz inequality for the positive quadrature weights yields
\[
\bigl|z^T\bigl(R_{\tilde f}(u)-R_{\tilde f}(v)\bigr)\bigr|
\le L\Bigl(\sum_q w_q\, z_h(x_q)^2\Bigr)^{\!1/2}
      \Bigl(\sum_q w_q\, (u_h-v_h)(x_q)^2\Bigr)^{\!1/2}.
\]
By \eqref{eq:quad-rep} the first factor equals $\norm{z}_M$ and the second $\norm{u-v}_M$, so
\[
\bigl|z^T\bigl(R_{\tilde f}(u)-R_{\tilde f}(v)\bigr)\bigr|
\le L\, \norm{z}_M\, \norm{u-v}_M.
\]
Taking the supremum over $z$ with $\norm{z}_M=1$ and using the $M$-norm identity
\[
\norm{w}_M=\sup_{\norm{z}_M=1}\ip{z}{w}_M
\]
with $w=M^{-1}(R_{\tilde f}(u)-R_{\tilde f}(v))$ gives \eqref{eq:LipM}.
\end{proof}

\begin{remark}
If $M$ and $R_{\tilde f}$ are assembled with two distinct positive quadrature rules, norm equivalence between the associated discrete $L^2$ norms introduces a constant $C_Q\ge 1$ and \eqref{eq:LipM} becomes $\norm{M^{-1}(R_{\tilde f}(u)-R_{\tilde f}(v))}_M \le C_Q L \norm{u-v}_M$. The start-up condition in Lemma~\ref{lem:startup} is then $S\ge C_Q L/2$. The sphere implementation used for the production runs assembles $M$ and $K$ with a rule of degree $2p+3$ and $R_{\tilde f}$, $\mathcal R_h$ with a rule of degree $4p+1$. On the production meshes, $C_Q$ in \eqref{eq:CQ} is at most $1+2.5\times10^{-6}$ (Section~\ref{sec:spatial}), and $S=600$ satisfies the start-up condition. The mass conservation of Section~\ref{sec:mass} and the exact balance \eqref{eq:balance} follow from the discrete update for any pair of rules, independently of $L$.
\end{remark}

\begin{lemma}[Discrete Taylor bound]\label{lem:Taylor}
Assume the common positive quadrature representation \eqref{eq:quad-rep}.
\begin{equation}
\mathcal R_h(u)-\mathcal R_h(v)
\le
(u-v)^T R_{\tilde f}(v)
+ \frac{L}{2}\norm{u-v}_M^2
\qquad \text{for all } u,v\in\R^{N_h}.
\label{eq:Taylor}
\end{equation}
\end{lemma}

\begin{proof}
For each quadrature point $x_q$, the standard one-dimensional $L$-smooth descent lemma for the $C^{1,1}$ function $\tilde F$ with $L$-Lipschitz derivative $\tilde f$ gives
\[
\tilde F\bigl(u_h(x_q)\bigr)-\tilde F\bigl(v_h(x_q)\bigr)
\le
\tilde f\bigl(v_h(x_q)\bigr)\bigl(u_h(x_q)-v_h(x_q)\bigr)
+ \frac{L}{2}\bigl(u_h(x_q)-v_h(x_q)\bigr)^2.
\]
Multiply by the positive weight $w_q$ and sum over $Q_h$. Using \eqref{eq:Rh} and the quadrature representation \eqref{eq:quad-rep} yields \eqref{eq:Taylor}.
\end{proof}

\begin{lemma}[Start-up energy decay]\label{lem:startup}
Under the hypotheses of Lemma~\ref{lem:Taylor} and the stabilization condition $S\ge L/2$ (with $L$ replaced by $C_QL$ for two rules, Remark after Lemma~\ref{lem:Lip}), every start-up step \eqref{eq:ie-1}--\eqref{eq:ie-2}, $k=1,\dots,n_0$, satisfies, with $d^k:=u^k-u^{k-1}$,
\begin{equation}
\mathcal E_h(u^k)-\mathcal E_h(u^{k-1})
+ \dt\, \norm{\mu^k}_K^2
+ \frac{\eps}{2}\norm{d^k}_K^2
+ \frac{S-L/2}{\eps}\norm{d^k}_M^2
\le 0.
\label{eq:startup-law}
\end{equation}
In particular,
\[
\mathcal E_h(u^{n_0})\le\dots\le\mathcal E_h(u^1)\le \mathcal E_h(u^0).
\]
More precisely, the exact balance
\begin{equation}
\mathcal E_h(u^k)-\mathcal E_h(u^{k-1})
+ \dt\, \norm{\mu^k}_K^2
+ \frac{\eps}{2}\norm{d^k}_K^2
+ \frac{S}{\eps}\norm{d^k}_M^2
= \frac{D_0^k}{\eps},
\qquad
D_0^k := \mathcal R_h(u^k)-\mathcal R_h(u^{k-1})-(d^k)^TR_{\tilde f}(u^{k-1}),
\label{eq:startup-balance}
\end{equation}
holds without any condition on $S$, and Lemma~\ref{lem:Taylor} gives $D_0^k\le \tfrac{L}{2}\norm{d^k}_M^2$.
\end{lemma}

\begin{proof}
Fix $k$ and set $d:=d^k$. Left-multiplying \eqref{eq:ie-1} by $\dt\,(\mu^k)^T$ gives
\begin{equation}
\ip{d}{\mu^k}_M + \dt\,\norm{\mu^k}_K^2 = 0.
\label{eq:startup-test1}
\end{equation}
Left-multiplying \eqref{eq:ie-2} by $d^T$ gives
\begin{equation}
\ip{d}{\mu^k}_M
=
\eps\,\ip{d}{u^k}_K
+ \frac{S}{\eps}\norm{d}_M^2
+ \frac1{\eps}\, d^T R_{\tilde f}(u^{k-1}).
\label{eq:startup-test2}
\end{equation}
Eliminate $\ip{d}{\mu^k}_M$ between \eqref{eq:startup-test1} and \eqref{eq:startup-test2}. The polarization identity gives
\[
\eps\,\ip{d}{u^k}_K
=
\frac{\eps}{2}\Bigl(
\norm{u^k}_K^2-\norm{u^{k-1}}_K^2+\norm{d}_K^2
\Bigr),
\]
and adding and subtracting $\eps^{-1}(\mathcal R_h(u^k)-\mathcal R_h(u^{k-1}))$ yields \eqref{eq:startup-balance}. By \eqref{eq:Taylor} with $u=u^k$ and $v=u^{k-1}$,
\[
\mathcal R_h(u^k)-\mathcal R_h(u^{k-1})
\le
d^T R_{\tilde f}(u^{k-1})
+ \frac{L}{2}\norm{d}_M^2,
\]
that is, $D_0^k\le\tfrac L2\norm{d}_M^2$. Substituting this bound into \eqref{eq:startup-balance} yields \eqref{eq:startup-law}.
\end{proof}

Let us next introduce notation to simplify the next stage of discussion.
Set
\[
\delta_0^n := u^n-u^{n-1},
\qquad
\delta_1^n := u^{n+1}-u^n,
\qquad
\delta_2^n := u^{n+1}-2u^n+u^{n-1},
\]
and
\[
B^n := \frac32 u^{n+1}-2u^n+\frac12 u^{n-1}.
\]
Then $u^{n+1}-u^{\mathrm{ext},n}=\delta_2^n$.

\begin{lemma}[BDF2 stabilization identity]\label{lem:stab-id}
For every symmetric positive semidefinite matrix $A\in\R^{N_h\times N_h}$ and every $n\ge n_0$,
\begin{equation}
\ip{B^n}{\delta_2^n}_A
= \frac12\bigl(\norm{\delta_1^n}_A^2-\norm{\delta_0^n}_A^2\bigr)
+ \norm{\delta_2^n}_A^2.
\label{eq:stab-id}
\end{equation}
\end{lemma}

\begin{proof}
Since $B^n=\delta_1^n+\tfrac12\delta_2^n$ and $\delta_2^n=\delta_1^n-\delta_0^n$,
\[
\ip{B^n}{\delta_2^n}_A
=\norm{\delta_1^n}_A^2-\ip{\delta_1^n}{\delta_0^n}_A
+\frac12\norm{\delta_2^n}_A^2.
\]
Using
$2\ip{\delta_1^n}{\delta_0^n}_A
=\norm{\delta_1^n}_A^2+\norm{\delta_0^n}_A^2-\norm{\delta_2^n}_A^2$
gives \eqref{eq:stab-id}.
\end{proof}

For every BDF2 step $n\ge n_0$ define the extrapolation defect
\begin{equation}
D^n := \mathcal R_h(u^{n+1})-\mathcal R_h(u^n)-(B^n)^T R_{\tilde f}^{\mathrm{ext},n},
\label{eq:Dn}
\end{equation}
the difference between the increment of the discrete potential and the work of the extrapolated reaction vector along $B^n$. This defect is the sole sign-indefinite term in the BDF2 energy balance.

\begin{assumption}[Nonlinear extrapolation defect bound]\label{ass:ext}
For the semi-implicit BDF2 sequence generated by \eqref{eq:bdf2-1}--\eqref{eq:bdf2-2}, assume there is a constant $\theta\ge0$, uniform in $n$ along that trajectory, such that for every $n\ge n_0$,
\begin{equation}
D^n
\le
(S-\theta)\norm{\delta_2^n}_M^2.
\label{eq:ext-defect}
\end{equation}
\end{assumption}

\begin{corollary}[Discrete chain rule for stabilized BDF2]\label{cor:bdf2-chain}
Under Assumption~\ref{ass:ext}, for every $n\ge n_0$,
\begin{align}
\frac{S}{\eps}\,\ip{B^n}{\delta_2^n}_M
+ \frac1{\eps} (B^n)^T R_{\tilde f}^{\mathrm{ext},n}
&\ge
\frac1{\eps}\Bigl(\mathcal R_h(u^{n+1})-\mathcal R_h(u^n)\Bigr)
\notag\\
&\quad+
\frac{S}{2\eps}\Bigl(\norm{\delta_1^n}_M^2-\norm{\delta_0^n}_M^2\Bigr)
+
\frac{\theta}{\eps}\norm{\delta_2^n}_M^2.
\label{eq:bdf2-ext-est}
\end{align}
\end{corollary}

\begin{proof}
Apply Lemma~\ref{lem:stab-id} with $A=M$ to obtain
\[
\frac{S}{\eps}\ip{B^n}{\delta_2^n}_M
=\frac{S}{2\eps}\bigl(\norm{\delta_1^n}_M^2-\norm{\delta_0^n}_M^2\bigr)
+\frac{S}{\eps}\norm{\delta_2^n}_M^2.
\]
Rearranging \eqref{eq:ext-defect}, dividing by $\eps$, and adding the result to the preceding identity gives \eqref{eq:bdf2-ext-est}.
\end{proof}

\begin{remark}
The stabilization contribution to \eqref{eq:bdf2-ext-est}, including the telescoping term
$(S/2\eps)(\norm{\delta_1^n}_M^2-\norm{\delta_0^n}_M^2)$ and the positive term
$(S/\eps)\norm{\delta_2^n}_M^2$, is the algebraic content of Lemma~\ref{lem:stab-id}. The scalar defect bound \eqref{eq:ext-defect} controls the discrepancy between the discrete chain rule and the explicit extrapolation $R_{\tilde f}^{\mathrm{ext},n}$. The combination of vector extrapolation $2f(u^n)-f(u^{n-1})$ and second-difference stabilization appears in \citep[Eq.~(6.22a)]{Shen2011}. Related modified-energy analyses use BDF2 with Douglas--Dupont regularization, state extrapolation, or variable steps~\citep{YanChenWangWise2018,WangYu2018,LiaoJiWangZhang2022}. Those estimates concern different nonlinear evaluations or stabilization terms. Equation~\eqref{eq:ext-defect} therefore remains an assumption for the surface finite element, quadrature-defined, truncated vector-extrapolated scheme. For the truncation \eqref{eq:trunc}, \eqref{eq:L2} gives $L=2$. Section~\ref{sec:defect} reports the trajectory-dependent defect quotient.
\end{remark}

For $n\ge1$ define the modified discrete Lyapunov functional
\begin{equation}
E_h^n :=
\frac{\eps}{4}\Bigl(\norm{u^n}_K^2 + \norm{2u^n-u^{n-1}}_K^2\Bigr)
+ \frac1{\eps}\mathcal R_h(u^n)
+ \frac{S}{2\eps}\norm{\delta_0^n}_M^2.
\label{eq:Lyap}
\end{equation}

\begin{corollary}[Exact BDF2 energy balance]\label{cor:balance}
Suppose that the sequence satisfies the BDF2 step \eqref{eq:bdf2-1}--\eqref{eq:bdf2-2} for $n\ge n_0$, and let $D^n$ be defined by \eqref{eq:Dn}. Then, for every $n\ge n_0$,
\begin{equation}
E_h^{n+1} - E_h^n
+ \dt\, \norm{\mu^{n+1}}_K^2
+ \frac{\eps}{4}\norm{\delta_2^n}_K^2
+ \frac{S}{\eps}\norm{\delta_2^n}_M^2
= \frac{D^n}{\eps}.
\label{eq:balance}
\end{equation}
\end{corollary}

\begin{proof}
Left-multiply \eqref{eq:bdf2-1} by $\dt\,(\mu^{n+1})^T$:
\begin{equation}
\ip{B^n}{\mu^{n+1}}_M + \dt\,\norm{\mu^{n+1}}_K^2 = 0.
\label{eq:test1}
\end{equation}
Left-multiply \eqref{eq:bdf2-2} by $(B^n)^T$:
\begin{equation}
\ip{B^n}{\mu^{n+1}}_M
= \eps\,\ip{B^n}{u^{n+1}}_K
+ \frac{S}{\eps}\,\ip{B^n}{\delta_2^n}_M
+ \frac1{\eps} (B^n)^T R_{\tilde f}^{\mathrm{ext},n}.
\label{eq:test2}
\end{equation}
Eliminate $\ip{B^n}{\mu^{n+1}}_M$ between \eqref{eq:test1} and \eqref{eq:test2}. Applying \eqref{eq:Gid} with $A=K$, $a=u^{n+1}$, $b=u^n$, and $d=u^{n-1}$ gives
\begin{align*}
\eps\,\ip{B^n}{u^{n+1}}_K
&=
\frac{\eps}{4}\Bigl(
\norm{u^{n+1}}_K^2 + \norm{2u^{n+1}-u^n}_K^2
- \norm{u^n}_K^2 - \norm{2u^n-u^{n-1}}_K^2
\Bigr)
\\
&\qquad+
\frac{\eps}{4}\norm{\delta_2^n}_K^2,
\end{align*}
Lemma~\ref{lem:stab-id} with $A=M$ gives
$\frac{S}{\eps}\ip{B^n}{\delta_2^n}_M
=\frac{S}{2\eps}\bigl(\norm{\delta_1^n}_M^2-\norm{\delta_0^n}_M^2\bigr)
+\frac{S}{\eps}\norm{\delta_2^n}_M^2$,
and by definition \eqref{eq:Dn},
$\eps^{-1}(B^n)^TR_{\tilde f}^{\mathrm{ext},n}
=\eps^{-1}\bigl(\mathcal R_h(u^{n+1})-\mathcal R_h(u^n)\bigr)-\eps^{-1}D^n$.
Collecting terms, and noting $\delta_1^n=\delta_0^{n+1}$, yields \eqref{eq:balance}.
\end{proof}

Every term on the left of \eqref{eq:balance} except the energy increment is nonnegative. The identity therefore shows exactly how much defect the scheme can absorb: $E_h^n$ decreases at step $n$ whenever $D^n\le S\norm{\delta_2^n}_M^2+\eps\dt\norm{\mu^{n+1}}_K^2+\tfrac{\eps^2}{4}\norm{\delta_2^n}_K^2$, which is weaker than \eqref{eq:ext-defect}. Section~\ref{sec:defect} reports both $D^n$ and the residual of \eqref{eq:balance} along computed trajectories.

\begin{theorem}\label{thm:lyap}
Under the hypotheses of Corollary~\ref{cor:balance} and Assumption~\ref{ass:ext}, for every $n\ge n_0$,
\begin{equation}
E_h^{n+1} - E_h^n
+ \dt\, \norm{\mu^{n+1}}_K^2
+ \frac{\eps}{4}\norm{\delta_2^n}_K^2
+ \frac{\theta}{\eps}\norm{\delta_2^n}_M^2
\le 0.
\label{eq:Lyap-law}
\end{equation}
In particular,
\[
E_h^{n+1} \le E_h^n
\qquad \text{for all } n\ge n_0.
\]
\end{theorem}

\begin{proof}
Insert \eqref{eq:ext-defect} into the right-hand side of \eqref{eq:balance}: $D^n/\eps\le \eps^{-1}(S-\theta)\norm{\delta_2^n}_M^2$. Moving this bound to the left cancels the term $(S/\eps)\norm{\delta_2^n}_M^2$ up to $(\theta/\eps)\norm{\delta_2^n}_M^2$, which is \eqref{eq:Lyap-law}. Equivalently, the estimate is Corollary~\ref{cor:bdf2-chain} applied inside the proof of Corollary~\ref{cor:balance}.
\end{proof}

\begin{remark}
Lemma~\ref{lem:startup} gives decay of the truncated discrete energy $\mathcal E_h$ in \eqref{eq:Eh} through the start-up steps. Under Assumption~\ref{ass:ext}, Theorem~\ref{thm:lyap} gives monotonicity of the modified discrete Lyapunov functional $E_h^n$ in \eqref{eq:Lyap} for $n\ge n_0$. The two energies differ by the quadratic stabilization term $(S/2\eps)\norm{\delta_0^n}_M^2$ and by the $G$-norm splitting in the gradient contribution. At the first BDF2 index $n_0$,
\[
E_h^{n_0}
=\mathcal E_h(u^{n_0})
+\frac{\eps}{4}\Bigl(\norm{2u^{n_0}-u^{n_0-1}}_K^2-\norm{u^{n_0}}_K^2\Bigr)
+\frac{S}{2\eps}\norm{u^{n_0}-u^{n_0-1}}_M^2.
\]
When the numerical solution stays in $[-1,1]$, the discrete potential $\mathcal R_h$ coincides with the quartic double-well energy. Once overshoots occur, $\mathcal E_h$ is the truncated energy advanced by the semi-implicit scheme. The numerical experiments report both $\mathcal E_h(u^n)$ and $E_h^n$.
\end{remark}

\section{Numerical Experiments}\label{sec:experiments}

The experiments verify accuracy, test the discrete structure laws, quantify
the limitations of production time steps, compare computational cost, and
examine the effect of geometry and topology. Table~\ref{tab:experiment-overview}
identifies the purpose and diagnostics of each group. Across the evolving
tests, the algebraic mass invariant is preserved to linear-solve accuracy.
The distinction from a quadrature-based mass diagnostic is stated below.
Quadrature order increases with polynomial degree, and the smooth
manufactured-solution test verifies exponential $p$-convergence. Energy
statements distinguish the truncated energy $\mathcal E_h$, the modified
functional $E_h^n$, and the physical energy of the separate viscous model.
In particular, the large-step sphere runs test conservative relaxation and
its accuracy limits, with monotonicity assessed from the recorded histories.

\begin{table}[!htbp]
\centering
\small
\setlength{\tabcolsep}{4pt}
\caption{Purpose of the numerical experiments. The common discretization
and solver settings are given in Section~\ref{sec:setup}.}
\label{tab:experiment-overview}
\begin{tabular}{>{\raggedright\arraybackslash}p{0.17\textwidth}>{\raggedright\arraybackslash}p{0.29\textwidth}>{\raggedright\arraybackslash}p{0.27\textwidth}>{\raggedright\arraybackslash}p{0.17\textwidth}}
\toprule
Purpose & Experiments & Diagnostics & Results \\
\midrule
Verify accuracy & Sphere eigenvalues, single-mode growth, manufactured solution
& Spectral errors, growth rates, temporal and spatial convergence
& Tables~\ref{tab:eigenvalues}--\ref{tab:h_convergence} \\[3pt]
Test discrete structure & Sphere mass and energy histories, extrapolation defect
& Mass drift, energy increments, defect quotient, balance residual
& Table~\ref{tab:energy_mono}, Figures~\ref{fig:structure}--\ref{fig:defect} \\[3pt]
Establish accuracy limits & Production-parameter growth, post-onset refinement, published sphere comparison
& Discrete growth rates, reference errors, time normalization
& Tables~\ref{tab:production-rates}--\ref{tab:zhang} \\[3pt]
Compare cost & SI/CS solvers and equal-dimension $hp$ configurations
& Stepping cost, residual contraction, assembly and solve costs
& Tables~\ref{tab:dt_comparison}--\ref{tab:hp_comparison}, Figures~\ref{fig:spinodal}--\ref{fig:hp_snapshots} \\[3pt]
Examine geometry and topology & Torus bands and coarsening, viscous torus, implicit genus-2 and genus-3 surfaces
& Interface energy, Fourier content, winding, unstable-mode counts
& Tables~\ref{tab:torus-bands}--\ref{tab:genus23}, Figures~\ref{fig:torus_plain}--\ref{fig:genus3_hero} \\
\bottomrule
\end{tabular}
\end{table}
\FloatBarrier

\subsection{Common discretization and solver settings}\label{sec:setup}

The experiments share the finite element spaces, quadrature rules, and
solver conventions below. Individual tests specify their parameter choices.

\subsubsection{Spaces, meshes, and quadrature}\label{sec:spatial}

The sphere mesh starts from the $20$ faces of an inscribed regular
icosahedron. Each refinement splits every triangle into four and projects
the new vertices radially onto $\Sphere$, giving $20\cdot4^r$ curvilinear
elements at level $r$. The $C^0$ hierarchic space uses vertex modes
$\lambda_i$, edge modes $\lambda_i\lambda_jP_{k-2}(\lambda_j-\lambda_i)$
for $k=2,\ldots,p$, and interior modes
\[
\varphi_{ij}=\lambda_1\lambda_2\lambda_3
P^{(2,2)}_{i-1}(2\lambda_2-1)P^{(2,2)}_{j-1}(2\lambda_3-1),
\qquad i,j\ge1,\quad i+j\le p-1.
\]
These span $\mathbb P_p$ on the reference triangle, of dimension
$(p+1)(p+2)/2$. Vertex coefficients are nodal values. Edge and interior
coefficients are bubble weights. Pointwise ranges are therefore evaluated
from $u_h$.

On the sphere, positive Xiao--Gimbutas rules of degrees $2p+3$ and $4p+1$
assemble $(K,M)$ and the reaction, potential, and load vectors, respectively.
The corresponding matrix-rule degrees are $2p+7$ for the torus coarsening
and viscous runs and $2p+10$ for the implicit surfaces. Their reaction rule
has degree $4p+1$. The single-mode validation, controlled torus bands, and
deterministic torus benchmark use a separate implementation with one
positive degree-$(2p+7)$ rule for all terms.

For two rules, let $M_R$ be the mass matrix assembled with the reaction
rule. The norm-equivalence constant in the remark after Lemma~\ref{lem:Lip} is
\begin{equation}
C_Q:=\max\Bigl\{1,\ \sup_{z\ne0}\frac{z^TM_Rz}{z^TMz}\Bigr\}
=\max\{1,\lambda_{\max}(M^{-1}M_R)\}.
\label{eq:CQ}
\end{equation}
Arnoldi iteration for the generalized pair $(M_R,M)$ gives
$C_Q-1=2.1\times10^{-7}$, $2.4\times10^{-6}$, and $2.4\times10^{-7}$ on
r4$p$3, r3$p$6, and r4$p$5, respectively. On r2$p$3 and r0$p$2 the values
are $4.6\times10^{-5}$ and $8.0\times10^{-3}$. The production values leave
the start-up condition $S\ge C_QL/2$ well satisfied at $S=600$.
For $L_0$ the constant-function load, the much smaller residual
$\norm{L_0-Mc}/\norm{L_0}$ ($1.9\times10^{-16}$ on r4$p$5 and $10^{-4}$
on r0$p$2) tests only the constant mode. The two mass quantities are related by
\[
m_Q^n:=L_0^Tu^n,\qquad m_M^n:=c^TMu^n,
\qquad m_Q^n-m_M^n=(L_0-Mc)^Tu^n.
\]
The invariant is $m_M^n$. The production tables report changes of $m_Q^n$.

We use $N_h=\dim V_h$ for scalar DOFs. The coupled system has $2N_h$
unknowns. The sphere counts satisfy
\begin{equation}
  N_h(r,p)=V_r+(p-1)E_r+\tfrac12(p-1)(p-2)F_r
          =10\cdot 4^{\,r}p^2+2,
\label{eq:dofcount}
\end{equation}
with $F_r=20\cdot4^r$, $E_r=3F_r/2$, and $V_r=F_r/2+2$.
Table~\ref{tab:meshes} contains three equal-dimension configurations with
the same $h/p$, and the r4$p$5 production mesh with $h/p=0.0099$ at
$\eps=0.01$. The former isolate the allocation of resolution between
$h$ and $p$. The latter resolves the nominal diffuse-interface width.

\begin{table}[t]
\centering
\caption{Mesh configurations: refinement level, polynomial order, element
count, global degrees of freedom, and the resolution scale $h/p$ with
$h=\sqrt{4\pi/N_{\mathrm{faces}}}$.  The first three configurations have
identical global dimension by \eqref{eq:dofcount} and are used for the
$hp$-refinement study of Section~\ref{sec:hp_study}. The last is the
production configuration for the $\eps=0.01$ experiments, chosen so that
$h/p$ resolves the diffuse interface.}\label{tab:meshes}
\begin{tabular}{ccccc}
\toprule
Refinement $r$ & $p$ & Faces & DOFs & $h/p$ \\
\midrule
4 & 3  & 5120 & 23\,042 & 0.0165 \\
3 & 6  & 1280 & 23\,042 & 0.0165 \\
2 & 12 & 320  & 23\,042 & 0.0165 \\
\midrule
4 & 5  & 5120 & 64\,002 & 0.0099 \\
\bottomrule
\end{tabular}
\end{table}

\subsubsection{Time integration and initial conditions}\label{sec:production-setup}

Every semi-implicit run uses vector reaction extrapolation
\eqref{eq:ext} and the truncation \eqref{eq:trunc}. The CS comparison uses
the untruncated cubic of Section~\ref{sec:cs}. Beyond $[-1,1]$ the two
solvers advance different discrete models. The viscous torus uses the
separate SAV modification specified in Section~\ref{sec:voigt-test}.

The sphere production runs use $\eps=0.01$, $S=600$, $T=2$, and two
start-up steps. Vertex values are drawn from $[-0.05,0.05]$ with seed
$42$, then the projected constant mode is subtracted to enforce zero mass.
The same initial field is projected separately for the equal-dimension
$hp$ comparison. The three production steps are
$\dt=5\times10^{-5}$, $2\times10^{-4}$, and $5\times10^{-4}$.
Single-mode and restarted accuracy tests state their alternative initial
conditions and final times explicitly. Other non-viscous tests use one start-up
step unless stated otherwise.

With two start-up steps, a $4\,000$-step run contains $3\,998$ BDF2
updates. A single-start-up rerun at $\dt=5\times10^{-4}$ gives the same
reported terminal energy $10.761$ and trajectory maximum $3.30$, with
$E_h^n$ decreasing at every BDF2 step. The manufactured test below checks
the asymptotic order for both start-up choices.

\subsubsection{Solvers and computational environment}\label{sec:environment}

Computations use the CPU of an Apple M3 Max workstation (16 cores,
128 GB unified memory, macOS). Wolfram Language 14.3 produced
Tables~\ref{tab:growth-rate}, \ref{tab:dt_comparison}, \ref{tab:hp_comparison},
\ref{tab:energy_mono}, and \ref{tab:torus-bands}--\ref{tab:genus23}. Version
15.0.1 produced Tables~\ref{tab:temporal}--\ref{tab:h_convergence},
\ref{tab:cs_residuals}, \ref{tab:production-rates}, and
\ref{tab:coarsening-refinement}, the balance diagnostics, and reruns of the
growth, band, and deterministic torus tests. Timing comparisons use the
version identified with each experiment.

The sparse $2N_h\times2N_h$ SI matrices for start-up and BDF2 are each
factored once by the built-in multifrontal $LU$ solver. The mass matrix is
factored by Cholesky. SI requires one back-substitution per step. CS reuses
the same factorizations in \eqref{eq:PSI}, stopping when the scaled
Euclidean block residual satisfies
\[
\bigl(\dt^{-2}\norm{R_1}^2+\norm{R_2}^2\bigr)^{1/2}<10^{-10},
\qquad R_1=\dt F_1,\quad R_2=-F_2,
\]
with at most $50$ corrections. Each correction assembles the cubic reaction
and performs one back-substitution. A final residual evaluation confirms
convergence. Stepping times exclude initial assembly, factorization, and
per-step diagnostics. The largest sphere system ($128\,004$ coupled
unknowns) and both implicit-surface systems fit within the 128 GB memory
budget. Typical solve costs are $21$ ms for r4$p$5, and $4.7$ and $8.3$ ms
for the genus-2 and genus-3 systems.

\FloatBarrier
\subsection{Verification of spatial and temporal accuracy}\label{sec:verification}

Sphere eigenvalues test spatial assembly, resolved modal growth tests
time integration, and a manufactured solution measures convergence of
both fields under $h$-, $p$-, and time-step refinement.

\subsubsection{Sphere eigenvalues}\label{sec:eigenvalues}

The generalized eigenproblem $K\phi=\lambda M\phi$ tests matrix assembly
against $\lambda_l=l(l+1)$ with multiplicity $2l+1$. On the $r=3$ mesh,
Table~\ref{tab:eigenvalues} shows spectral convergence in $p$, with the
first four distinct eigenvalues at floating-point accuracy by $p=6$.
The multiplicities $1,3,5,7$ are recovered by shift-invert subspace
iteration with a sparse factorization of $K-\sigma M$.

\begin{table}[t]
\centering
\caption{Eigenvalue errors $|\lambda_h - \lambda_{\mathrm{exact}}|$ for the first four distinct eigenvalues of $-\LapS$ on the $r=3$ mesh (1280 faces).  Multiplicity is recovered correctly in all cases.}\label{tab:eigenvalues}
\begin{tabular}{ccccc}
\toprule
$p$ & $\lambda_1 = 0$ & $\lambda_{2\text{--}4} = 2$ & $\lambda_{5\text{--}9} = 6$ & $\lambda_{10\text{--}16} = 12$ \\
\midrule
2  & $1.1\times10^{-15}$ & $9.2\times10^{-6}$  & $1.1\times10^{-4}$  & $7.5\times10^{-4}$  \\
4  & $2.3\times10^{-15}$ & $2.1\times10^{-11}$ & $4.5\times10^{-10}$ & $5.1\times10^{-9}$  \\
6  & $6.3\times10^{-15}$ & $1.4\times10^{-14}$ & $2.5\times10^{-14}$ & $5.3\times10^{-14}$ \\
\bottomrule
\end{tabular}
\end{table}

\subsubsection{Resolved linear growth}\label{sec:spectral-growth}

For a constant background $\bar u$, the linearized growth rate of a
Laplace--Beltrami mode is
\begin{equation}
\sigma(\lambda)
=\lambda\left(\frac{1-3\bar u^2}{\eps}-\eps\lambda\right).
\label{eq:growth-rate}
\end{equation}
and, when $3\bar u^2<1$, its unstable band and maximum are
\begin{equation}
0<\lambda<\frac{1-3\bar u^2}{\eps^2},
\qquad
\lambda_*=\frac{1-3\bar u^2}{2\eps^2},
\qquad
\sigma_*=\frac{(1-3\bar u^2)^2}{4\eps^3}.
\label{eq:spinodal-scales}
\end{equation}
The associated scales are $\ell_*=2\pi/\sqrt{\lambda_*}$ and
$\tau_*=1/\sigma_*$ \citep{MarenduzzoOrlandini2013}.

To test temporal growth independently of production coarsening, set
$\eps=1/\sqrt{24}$, $S=2$, $r=2$, $p=4$, $N_h=2562$, $T=0.1$, and
$u_h^0=10^{-3}P_hY_{3m}$. The fastest eigenspace then has $\lambda_*=12$.
Table~\ref{tab:growth-rate} reports the least-squares slope of the
logarithmic modal amplitude on $t\in[0.02,0.08]$ for $m=0$. Its error
decreases with each time-step halving. At $\dt=2.5\times10^{-4}$ the
relative spread across $m=-2,0,1,3$ is $9.9\times10^{-7}$. The maximum
mass drift over the three steps is $7.8\times10^{-17}$.

\begin{table}[pos=!htbp]
\centering
\caption{Single-mode sphere growth test. The analytical rate is
$\sigma(12)=12\sqrt6=29.3938769$.}
\label{tab:growth-rate}
\begin{tabular}{ccc}
\toprule
$\dt$ & measured $\sigma_h$ & relative error \\
\midrule
$5.00\times10^{-4}$ & 29.3582341 & $1.21\times10^{-3}$ \\
$2.50\times10^{-4}$ & 29.3846735 & $3.13\times10^{-4}$ \\
$1.25\times10^{-4}$ & 29.3913557 & $8.58\times10^{-5}$ \\
\bottomrule
\end{tabular}
\end{table}

\subsubsection{Manufactured solution}\label{sec:manufactured}

The forced problem tests convergence of both fields with $S=12$,
$\eps=0.1$, one start-up step, and
\[
u_{\mathrm{ex}}(\mathbf x,t)=\bar u+a(t)x_1,\qquad
a(t)=0.1e^{-t},\qquad \bar u=0.7,
\qquad
\mu_{\mathrm{ex}}=2\eps a(t)x_1+\eps^{-1}(u_{\mathrm{ex}}^3-u_{\mathrm{ex}}).
\]
Here $u_{\mathrm{ex}}\in[0.6,0.8]$, so truncation is inactive and
$f'(u_{\mathrm{ex}})\in[0.08,0.92]$. Nonconstant modes of the
constant-background linearization decay since $3\bar u^2>1$.
With $\alpha=1$, direct substitution using
$\LapS x_1=-2x_1$, $\LapS x_1^2=2-6x_1^2$, and
$\LapS x_1^3=6x_1-12x_1^3$ gives
\begin{equation}
\partial_t u - \LapS \mu = g_u,
\qquad
\mu + \eps \LapS u - \eps^{-1}(u^3-u) = 0,
\label{eq:forced}
\end{equation}
\begin{equation}
g_u
=-\frac{6a^2\bar u}{\eps}
+\Bigl[\Bigl(4\eps-\frac{2}{\eps}-\alpha+\frac{6\bar u^2}{\eps}\Bigr)a-\frac{6a^3}{\eps}\Bigr]x_1
+\frac{18a^2\bar u}{\eps}\,x_1^2
+\frac{12a^3}{\eps}\,x_1^3,
\label{eq:gu}
\end{equation}
The forcing is independent of the numerical states and stencil. Its load
$F_i(t)=\sum_qw_qg_u(x_q,t)\varphi_i(x_q)$ uses the reaction quadrature
and is inserted at $t^{n+1}$ into the first block row. The initial datum is
the $L^2$ projection. Final errors are measured as
$[\sum_qw_q(u_h-u_{\mathrm{ex}})^2(x_q)]^{1/2}$, and similarly for $\mu$.

Table~\ref{tab:temporal} fixes $r=3$, $p=4$, and $T=0.1$. After the first
refinement, rates are $1.95$--$2.09$ for $u$ and $1.96$--$2.01$ for $\mu$.
Two start-up steps give rates $1.9$--$2.0$ over the last three refinements,
with $u$ errors larger by factors $1.2$--$1.6$, and hence retain second order.
The stabilization contributes an $O(S\dt^2/\eps)$ consistency term.

\begin{table}[t]
\centering
\caption{Temporal convergence of the manufactured solution with the analytic forcing \eqref{eq:gu} ($r=3$, $p=4$, $\eps=0.1$, $S=12$, $T=0.1$, one start-up step).}\label{tab:temporal}
\begin{tabular}{ccccc}
\toprule
$\dt$ & $\|u - u_h\|_{L^2}$ & rate & $\|\mu - \mu_h\|_{L^2}$ & rate \\
\midrule
$2\times10^{-2}$ & $9.23\times10^{-3}$ & -- & $2.33\times10^{-2}$ & -- \\
$1\times10^{-2}$ & $3.09\times10^{-3}$ & 1.58 & $1.29\times10^{-2}$ & 0.85 \\
$5\times10^{-3}$ & $7.26\times10^{-4}$ & 2.09 & $3.32\times10^{-3}$ & 1.96 \\
$2\times10^{-3}$ & $1.13\times10^{-4}$ & 2.03 & $5.25\times10^{-4}$ & 2.01 \\
$1\times10^{-3}$ & $2.93\times10^{-5}$ & 1.95 & $1.35\times10^{-4}$ & 1.96 \\
\bottomrule
\end{tabular}
\end{table}

The spatial tests use $\dt=5\times10^{-6}$ and $T=0.05$, with temporal
error about $10^{-9}$. At fixed $r=2$, increasing $p$ reduces the $u$
error by factors $9$--$27$ per degree (Table~\ref{tab:p_convergence}).
At fixed $p=1,2,3$, the finest $h$-refinement rates approach $p+1$ for
both fields (Table~\ref{tab:h_convergence}).

\begin{table}[t]
\centering
\caption{Spatial convergence under $p$-refinement ($r=2$, $\dt = 5\times10^{-6}$, $T=0.05$).}\label{tab:p_convergence}
\begin{tabular}{cccc}
\toprule
$p$ & DOFs & $\|u - u_h\|_{L^2}$ & $\|\mu - \mu_h\|_{L^2}$ \\
\midrule
1 & 162 & $1.07\times10^{-3}$ & $6.79\times10^{-3}$ \\
2 & 642 & $9.56\times10^{-5}$ & $5.82\times10^{-4}$ \\
3 & 1\,442 & $3.60\times10^{-6}$ & $2.66\times10^{-5}$ \\
4 & 2\,562 & $3.14\times10^{-7}$ & $2.16\times10^{-6}$ \\
5 & 4\,002 & $1.64\times10^{-8}$ & $1.29\times10^{-7}$ \\
6 & 5\,762 & $1.78\times10^{-9}$ & $1.16\times10^{-8}$ \\
\bottomrule
\end{tabular}
\end{table}

\begin{table}[t]
\centering
\caption{Spatial convergence under $h$-refinement ($\dt = 5\times10^{-6}$, $T=0.05$).  The characteristic mesh size is $h \sim \sqrt{4\pi/N_{\mathrm{faces}}}$.}\label{tab:h_convergence}
\begin{tabular}{cccccc}
\toprule
& & \multicolumn{2}{c}{$\|u - u_h\|_{L^2}$} & \multicolumn{2}{c}{$\|\mu - \mu_h\|_{L^2}$} \\
\cmidrule(lr){3-4} \cmidrule(lr){5-6}
$p$ & Level ($r$) & error & rate & error & rate \\
\midrule
\multirow{4}{*}{1}
 & 0 & $1.45\times10^{-2}$ & -- & $1.05\times10^{-1}$ & -- \\
 & 1 & $4.18\times10^{-3}$ & 1.80 & $2.69\times10^{-2}$ & 1.96 \\
 & 2 & $1.07\times10^{-3}$ & 1.96 & $6.79\times10^{-3}$ & 1.99 \\
 & 3 & $2.70\times10^{-4}$ & 1.99 & $1.71\times10^{-3}$ & 1.99 \\
\midrule
\multirow{4}{*}{2}
 & 0 & $5.10\times10^{-3}$ & -- & $2.98\times10^{-2}$ & -- \\
 & 1 & $7.46\times10^{-4}$ & 2.77 & $4.48\times10^{-3}$ & 2.73 \\
 & 2 & $9.56\times10^{-5}$ & 2.96 & $5.82\times10^{-4}$ & 2.95 \\
 & 3 & $1.20\times10^{-5}$ & 2.99 & $7.35\times10^{-5}$ & 2.99 \\
\midrule
\multirow{4}{*}{3}
 & 0 & $8.43\times10^{-4}$ & -- & $5.88\times10^{-3}$ & -- \\
 & 1 & $5.69\times10^{-5}$ & 3.89 & $4.16\times10^{-4}$ & 3.82 \\
 & 2 & $3.60\times10^{-6}$ & 3.98 & $2.66\times10^{-5}$ & 3.97 \\
 & 3 & $2.25\times10^{-7}$ & 4.00 & $1.67\times10^{-6}$ & 3.99 \\
\bottomrule
\end{tabular}
\end{table}

The forcing has zero discrete mean: the icosahedral rule gives
$c^TL_1,c^TL_3=0$ and $c^TL_2=c^TL_0/3$ to $2\times10^{-15}$ for the
loads of $x_1^k$. The recorded $|c^TF(t^n)|$ is at most
$4\times10^{-15}$, and $m_M^n$ changes by at most $1.1\times10^{-11}$.
The reaction-rule diagnostic $m_Q^n$ has a distinct initial offset and
subsequent drift on coarse meshes: on r0$p$2 these are $8.0\times10^{-4}$
and $3.2\times10^{-6}$, respectively. Maximum changes are
$1.5\times10^{-9}$ on r1$p$2, $9.6\times10^{-8}$ on r0$p$3,
$2.1\times10^{-11}$ on r1$p$3, and at most $2.1\times10^{-12}$ elsewhere.
The solution-error norms above are evaluated directly with the reaction rule.

\FloatBarrier
\subsection{Tests of discrete conservation and energy laws}\label{sec:structure-tests}

Mass and energy histories test the discrete invariants at production
parameters. The extrapolation defect and balance residual connect the
observed energy changes to the sufficient condition in
Theorem~\ref{thm:lyap} and the identity in Corollary~\ref{cor:balance}.

\subsubsection{Mass and energy histories}\label{sec:structure}

\begin{figure}
\centering
\includegraphics[width=\textwidth]{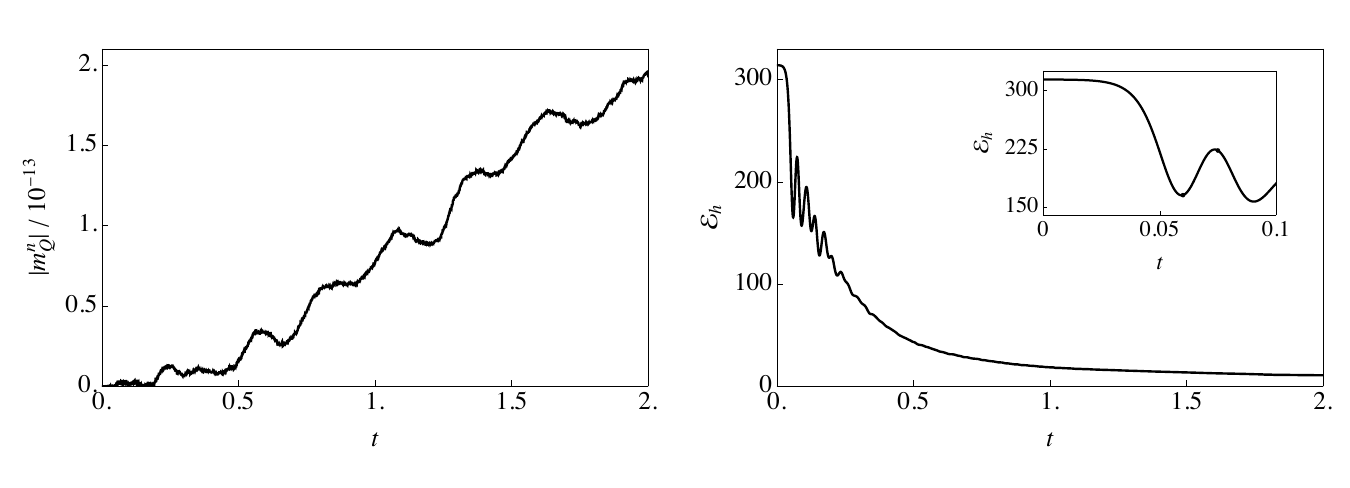}
\caption{\CapStructure}\label{fig:structure}
\end{figure}

The r4$p$5 production runs test the mass invariant and distinguish the
two energy diagnostics. Figure~\ref{fig:structure} shows the quadrature
mass and truncated energy at every step of the version-15.0.1 rerun at
$\dt=5\times10^{-4}$. The absolute mass remains at $O(10^{-13})$ for
the zero-mass initial field.
Table~\ref{tab:energy_mono} counts energy increments at every step of
the original version-14.3 runs for all three production steps.

\begin{table}[t]
\centering
\caption{Every-step monitoring of the two energy functionals for the
semi-implicit scheme in the original Wolfram Language~14.3 runs
($r=4$, $p=5$, $\eps=0.01$, $S=600$, $T=2$). $E_h^n$ is
the modified Lyapunov functional \eqref{eq:Lyap}. $\mathcal E_h$ is the
truncated discrete energy \eqref{eq:Eh}.  ``Increases'' counts steps with a
positive increment.}\label{tab:energy_mono}
\begin{tabular}{ccccccc}
\toprule
& & \multicolumn{2}{c}{$E_h^n$} & \multicolumn{2}{c}{$\mathcal E_h$} \\
\cmidrule(lr){3-4}\cmidrule(lr){5-6}
$\dt$ & Steps & increases & max increment & increases & max increment \\
\midrule
$5\times10^{-4}$ & 4\,000  & $0$    & --                 & 339 \ (8.5\,\%)  & $3.25$ \\
$2\times10^{-4}$ & 10\,000 & $0$    & --                 & 893 \ (8.9\,\%)  & $3.19$ \\
$5\times10^{-5}$ & 40\,000 & 4\,361 & $1.1\times10^{-7}$  & 9\,310 \ (23.3\,\%) & $3.01$ \\
\bottomrule
\end{tabular}
\end{table}

The modified energy $E_h^n$ decreases at every step at the two larger
steps. At $\dt=5\times10^{-5}$ its largest positive increment is
$1.1\times10^{-7}$, about $10^{-8}$ relative, confined to the terminal
stage $t\gtrsim1.47$ where the block factorization reports ill-conditioning.
The version-15.0.1 rerun has only $129$ positive increments, all below
$5.3\times10^{-10}$ at $t>1.957$. Small-increment counts depend on
linear-solve accuracy. Its $\mathcal E_h$ counts are $344$, $874$, and
$8\,779$ at decreasing step sizes, with the maximal $\mathcal E_h$
increments unchanged to three digits.

Decay of $E_h^n$ does not imply stepwise decay of $\mathcal E_h$.
The truncated energy $\mathcal E_h$ has substantial positive increments.
At $\dt=5\times10^{-4}$ it rises from $165$ at $t=0.06$ to $223$ at
$t=0.075$, before further oscillations and relaxation to $10.76$.
The largest single increment is $3.25$ at $t=0.067$, where $q^n\approx16$,
after the defect-bound violation described below. The stabilization and
$G$-norm terms separate $E_h^n$ from $\mathcal E_h$ by $39$ at $t=0.3$
and $16$ at $t=0.5$.

Overshoots activate the truncated model on a substantial part of the
sphere. At decreasing production steps, peak active-area fractions are
$74\%$, $74\%$, and $73\%$, and terminal fractions are $47\%$, about
$45\%$, and $46\%$. Trajectory maxima are $3.30$, $3.17$, and $3.47$,
despite terminal ranges close to $[-1,1]$ at the smaller steps
(Table~\ref{tab:dt_comparison}).

\subsubsection{Extrapolation defect and exact balance}\label{sec:defect}

\begin{figure}
\centering
\begin{subfigure}[t]{0.49\textwidth}
\centering
\includegraphics[width=\textwidth]{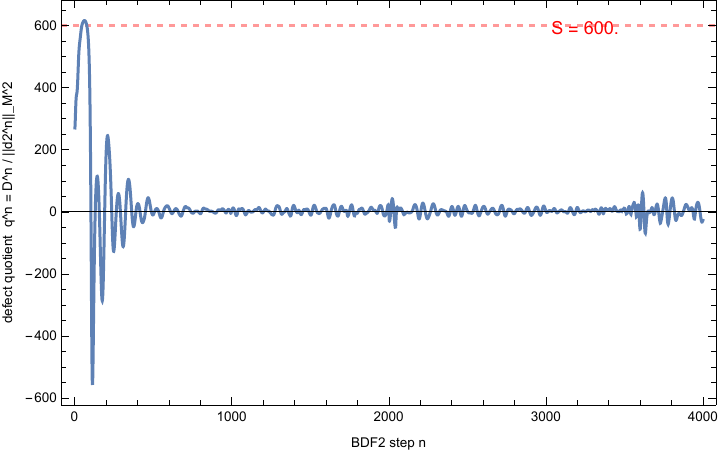}
\caption{Defect quotient $q^n$ of \eqref{eq:qn} (solid) against the
stabilization parameter $S=600$ (dashed). At each plotted step, $q^n\le S$
is equivalent to a nonnegative stepwise margin $\theta^n=S-q^n$. The bound
holds at $99.2\%$ of steps. The quotient exceeds $S$ by at most
$\approx2.6\%$ in the brief onset window $t\in[0.024,0.039]$ and is well
below $S$ thereafter.}
\label{fig:defect-q}
\end{subfigure}
\hfill
\begin{subfigure}[t]{0.49\textwidth}
\centering
\includegraphics[width=\textwidth]{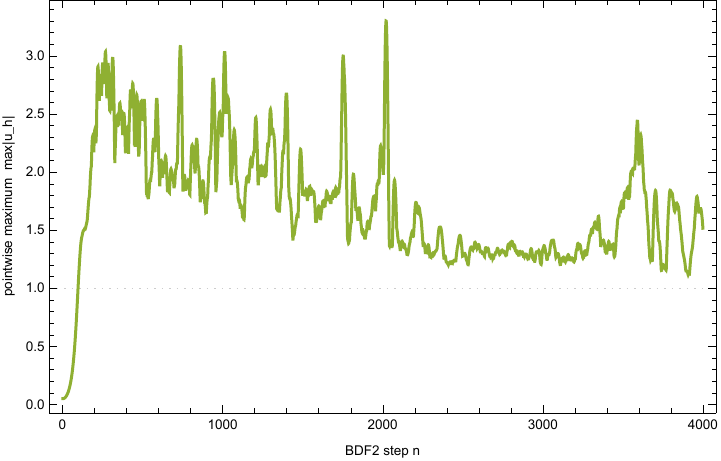}
\caption{Pointwise maximum $\max|u_h|$ (evaluated at the quadrature nodes)
at the same steps.  At the $q^n$ peak it is $\approx0.31$, so the
truncation is inactive in the window where the bound is violated.  The field
later leaves $[-1,1]$, reaching $3.30$ at $t\approx1.01$ before relaxing to
$1.51$ at $T=2$ (Table~\ref{tab:dt_comparison}). The truncation is active
there, but the defect bound holds throughout that stage.}
\label{fig:defect-overshoot}
\end{subfigure}
\caption{Empirical test of the extrapolation-defect bound,
Assumption~\ref{ass:ext}, for the semi-implicit BDF2 run at $\eps=0.01$,
$\dt=5\times10^{-4}$, $S=600$.}\label{fig:defect}
\end{figure}

For each step with $\norm{\delta_2^n}_M>0$, define the defect quotient
\begin{equation}
q^n := \frac{D^n}{\norm{\delta_2^n}_M^2}.
\label{eq:qn}
\end{equation}
For steps with $\delta_2^n=0$, check the condition $D^n\le0$ in
Assumption~\ref{ass:ext} directly. If there are steps with $\delta_2^n\ne0$,
the assumption also requires $\sup_{n:\,\delta_2^n\ne0}q^n\le S$, and one
can take $\theta=S-\sup_{n:\,\delta_2^n\ne0}q^n$.
At $\dt=5\times10^{-4}$ the stepwise bound holds at $3967$ of $3998$ BDF2
steps ($99.2\%$). The $31$ exceptions lie in $t\in[0.024,0.039]$, with
peak $q^n=615.5$ and $\max|u_h|\approx0.31$: truncation is inactive
during this window. The later overshoot reaches $3.30$ at $t\approx1.01$
while $q^n<S$ (Figure~\ref{fig:defect}).

The exact balance \eqref{eq:balance} closes to $4.5\times10^{-13}$ over
the BDF2 phase and the start-up balance to $3.1\times10^{-14}$, against
energies between $10$ and $314$. At $\dt=2\times10^{-4}$ and
$5\times10^{-5}$ the BDF2 residual bounds are $4.5\times10^{-13}$ and
$5.5\times10^{-13}$. Onset violations last $34$ and $42$ steps, with
peak quotients $619$ and $633$. Dissipation exceeds the excess defect in
the three onset windows by factors of at least $5.1$, $8.2$, and $13.8$,
explaining the observed decrease of $E_h^n$ despite failure of the
stronger sufficient bound.

At the smallest step, another $1\,288$ quotient exceedances occur for
$t>1.5$, with $\norm{\delta_2^n}_M^2\approx10^{-14}$ and excess defect
below $4\times10^{-9}$. Roundoff dominates the quotient in this regime.
The onset peaks apply to $S=600$. Assessing another stabilizer requires
recomputing the trajectory and its defect.

\FloatBarrier
\subsection{Accuracy limits of the production time steps}\label{sec:production-accuracy}

Linear growth rates quantify the effect of stabilization during onset.
A common-state restart measures time-step sensitivity during coarsening,
and a time normalization places these runs alongside published sphere
computations.

\subsubsection{Linear growth at production parameters}

For $\bar u=0$ and $\eps=0.01$, \eqref{eq:spinodal-scales} gives
$\lambda_*=5000$, $\ell_*=0.0889$, and $\tau_*=4\times10^{-6}$.
Sphere modes through $l=99$ are unstable, with the fastest near $l=70$.
The production ratios $\dt/\tau_*=12.5,50,125$ motivate a separate
test of the discrete growth rates. For $K\phi=\lambda M\phi$, the
linearized BDF2 amplification polynomial is
\begin{equation}
\Bigl(\tfrac{3}{2\dt}+\eps\lambda^2+\tfrac{S}{\eps}\lambda\Bigr)z^2
-\Bigl(\tfrac{2}{\dt}+\tfrac{2(S+1)}{\eps}\lambda\Bigr)z
+\Bigl(\tfrac{1}{2\dt}+\tfrac{S+1}{\eps}\lambda\Bigr)=0,
\label{eq:bdf2-amplification}
\end{equation}
Its dominant root defines $\sigma_{\dt}=\log z_1/\dt$.
Table~\ref{tab:production-rates} uses r4$p$5, $S=600$, and
$u_h^0=10^{-4}P_hY_{l0}$ for $l=20,70$. Rates are fitted to
$\log|a_n|$, $a_n=\int_{\Sphere}u_h^nY_{l0}\,\dd S$, over the second
half of a $200$-step run, after the parasitic root
$|z_2|\approx0.97$ has decayed and while $|u_h|\le0.3$.

\begin{table}[t]
\centering
\caption{Linear growth at production parameters ($\eps=0.01$, $S=600$, $r=4$, $p=5$). $\sigma$ is the continuous rate \eqref{eq:growth-rate}, $\sigma_{\dt}$ the dominant root of \eqref{eq:bdf2-amplification}, and $\sigma_h$ the rate measured in the semi-implicit computation.}\label{tab:production-rates}
\begin{tabular}{cccccc}
\toprule
$l$ & $\lambda$ & $\dt$ & $\sigma$ & $\sigma_{\dt}$ & $\sigma_h$ \\
\midrule
20 & 420 & $5\times10^{-4}$ & 40\,236 & 81.40 & 81.39 \\
20 & 420 & $2\times10^{-4}$ & 40\,236 & 203.2 & 203.2 \\
20 & 420 & $5\times10^{-5}$ & 40\,236 & 806.5 & 806.4 \\
70 & 4970 & $5\times10^{-4}$ & 249\,991 & 58.71 & 58.65 \\
70 & 4970 & $2\times10^{-4}$ & 249\,991 & 146.7 & 146.6 \\
70 & 4970 & $5\times10^{-5}$ & 249\,991 & 586.4 & 585.9 \\
70 & 4970 & $1\times10^{-5}$ & 249\,991 & 2\,918 & 2\,916 \\
70 & 4970 & $1\times10^{-6}$ & 249\,991 & 27\,662 & 27\,639 \\
\bottomrule
\end{tabular}
\end{table}

Measured and discrete rates agree within $0.1\%$, but differ greatly
from the continuous rates. At $\dt=5\times10^{-4}$, unstable modes
$1\le l\le70$ grow at rates $59$--$82$, falling to $8.2$ at $l=99$.
Stabilization flattens the dispersion relation and makes growth per step
nearly independent of $\dt$ across the tested production steps. Accordingly,
the onset energy fall and rebound occur over the same $30$ steps at
$\dt=5\times10^{-4}$ and $2\times10^{-4}$, with their time axes scaled
by the step ratio. Recovering the fastest continuous rate within $10\%$
requires about $\dt\le10^{-8}$ at $S=600$ ($10^{-7}$ at $S=2$), or
$5\times10^3$--$5\times10^4$ more steps than the production runs at
$S=600$. Their onset kinetics are therefore governed by the discretization.

\subsubsection{Post-onset time-step refinement}

All rows of Table~\ref{tab:coarsening-refinement} start from the common
SI production state $u_h(0.1)$ computed at $\dt=5\times10^{-4}$ and
advance to $T=0.5$, restarting with two backward-Euler steps of their own
size. The SI reference uses $\dt=5\times10^{-6}$ ($80\,000$ steps).
All errors are relative to this numerical reference, whose own
discretization error remains unquantified.

\begin{table}[t]
\centering
\footnotesize
\setlength{\tabcolsep}{3.5pt}
\caption{Post-onset refinement from the common state $u_h(0.1)$ to $T=0.5$ ($r=4$, $p=5$, $\eps=0.01$, $S=600$). Errors are against the semi-implicit run at $\dt=5\times10^{-6}$; $\norm{\cdot}_M$ is the discrete $L^2(\Sphere)$ norm, the relative error divides by $\norm{u_{\mathrm{ref}}(T)}_M=3.50$, and $\Delta\mathcal E_h=\mathcal E_h(T)-\mathcal E_h^{\mathrm{ref}}(T)$. ``corr./step'' is the mean number of convex-splitting corrections (back-substitutions) per step. Wall times are stepping costs: the per-step cost measured on an otherwise idle machine (Section~\ref{sec:environment}, Wolfram Language~15.0.1; $36$\,ms per semi-implicit step at every $\dt$, $250$\,ms per convex-splitting step at $5.6$ corrections, scaled in proportion to the number of residual evaluations) times the number of steps.}\label{tab:coarsening-refinement}
\begin{tabular}{lcccccccc}
\toprule
Solver & $\dt$ & Steps & $\norm{u_h-u_{\mathrm{ref}}}_M$ & rel. & $\mathcal E_h(T)$ & $\Delta\mathcal E_h$ & corr./step & Wall (min) \\
\midrule
SI & $5\times10^{-4}$ & 800 & $3.78$ & 108.1\% & 39.645 & $+29.34$ & -- & 0.48 \\
SI & $2\times10^{-4}$ & 2\,000 & $2.78$ & 79.5\% & 17.010 & $+6.70$ & -- & 1.2 \\
SI & $5\times10^{-5}$ & 8\,000 & $0.526$ & 15.0\% & 10.343 & $+0.037$ & -- & 4.8 \\
SI & $2\times10^{-5}$ & 20\,000 & $0.214$ & 6.1\% & 10.304 & $-0.002$ & -- & 12.0 \\
SI (reference) & $5\times10^{-6}$ & 80\,000 & -- & -- & 10.306 & -- & -- & 48.0 \\
CS & $5\times10^{-4}$ & 800 & $3.74$ & 106.8\% & 37.394 & $+27.09$ & 5.6 & 3.3 \\
CS & $2\times10^{-4}$ & 2\,000 & $2.78$ & 79.5\% & 16.716 & $+6.41$ & 5.4 & 8.0 \\
CS & $5\times10^{-5}$ & 8\,000 & $0.562$ & 16.1\% & 10.331 & $+0.025$ & 4.6 & 28.0 \\
\bottomrule
\end{tabular}
\end{table}

SI relative errors decrease from $108\%$ and $80\%$ at the two larger
production steps to $15\%$ at $5\times10^{-5}$ and $6.1\%$ at
$2\times10^{-5}$, relative to this reference. The smaller two errors
exhibit a pre-asymptotic first-order trend. The start-up choice retains
asymptotic second order (Table~\ref{tab:temporal}). Pointwise errors of
$1.2$--$2.5$ are dominated by displaced interfaces. At the larger steps,
restart overshoots and delayed coarsening also leave substantially higher
terminal energies.

CS gives similar distances from the SI reference at equal steps, at six
to seven times the cost. These distances measure each scheme against the
same SI reference. Comparing the SI and CS states directly requires a
pairwise error, since their discrete models differ beyond $[-1,1]$.
The order-one errors refer to this restarted experiment at $T=0.5$.
The uninterrupted $T=2$ production runs provide qualitative phase-separated
configurations and tests of discrete conservation and modified-energy behavior.
Their configurations at the final simulated time depend on the time step.
The present comparisons establish neither time-accurate continuum kinetics
nor independently validated morphological selection at the production steps.

\subsubsection{Relation to published sphere computations}\label{sec:zhang}

Table~\ref{tab:zhang} places the production regime alongside
\citet{ZhangWangRong2011}. Their Eq.~(3.27) uses the same stabilized
SI-BDF2 integrator, with prolate spectral elements and the untruncated
cubic. Their equation $\partial_tu=\Delta(-\eps^2\Delta u+f(u))$
implies $t_{\mathrm{present}}=\eps t_{\mathrm{ZWR}}$: their step and
final time become $10^{-5}$ and $0.02$ in the present units. The comparison
is qualitative, since spatial spaces, reactions, durations, and steps differ.

\begin{table}[t]
\centering
\caption{Reference comparison with the prolate-element computation of \citet{ZhangWangRong2011} at $\eps=0.01$. Times and time steps are converted to the present normalization by $t_{\mathrm{present}}=\eps\,t_{\mathrm{ZWR}}$. Table~\ref{tab:energy_mono} reports the energy increments for the present runs.}\label{tab:zhang}
\begin{tabular}{lcc}
\toprule
& \citet{ZhangWangRong2011} & Present work (r4$p$5, SI) \\
\midrule
Spatial discretization & prolate spectral elements, $N=32$ & mixed $hp$-FEM \\
Unknowns & ${\approx}\,6\,144$ points & $64\,002$ DOFs \\
Time stepping & stabilized SI-BDF2, Eq.~(3.27) & stabilized SI-BDF2, \eqref{eq:bdf2-1}--\eqref{eq:bdf2-2} \\
Reaction & $u^3-u$ & truncated $\tilde f$ of \eqref{eq:trunc} \\
$\dt$ (present units) & $10^{-5}$ & $5\times10^{-5}$ to $5\times10^{-4}$ \\
Final time (present units) & $0.02$ & $2$ \\
Mass diagnostic & $\max_k|\!\int u^k|=2.94\times10^{-8}$ & $|\delta m|_{\max}=1.8\times10^{-12}$ \\
Energy diagnostic & decreasing (their Fig.~10)
& \makecell{Stepwise $E_h^n$ decay at the two larger steps\\
  Small late increases at the smallest step} \\
\bottomrule
\end{tabular}
\end{table}

The linear second-order scheme of \citet{LiKimWang2017} provides an
alternative using Crank--Nicolson and stabilized splitting, with an energy
law under a bounded-state hypothesis and a sufficiently large stabilizer.
Transferring that law to the present $hp$ space and quadrature would
require a separate argument. Here the exact BDF2 balance and the
trajectory-dependent sufficient condition are tested directly.

\FloatBarrier
\subsection{Computational cost and resolution allocation}\label{sec:cost}

The solver comparison measures the cost of the nonlinear corrections.
The equal-dimension $hp$ comparison separates reaction assembly and
linear solves to quantify the cost of allocating resolution between
mesh refinement and polynomial degree.

\subsubsection{Semi-implicit and convex-splitting solvers}\label{sec:solver_comparison}

Table~\ref{tab:dt_comparison} compares SI and CS on the same r4$p$5
mesh, initial field, quadrature choices, and hardware, with the settings of
Section~\ref{sec:setup}. CS/SI stepping-time ratios are $4.2$, $5.1$,
and $5.8$ as $\dt$ increases. These ratios compare implementation costs at
fixed mesh and time step. The potentials differ outside $[-1,1]$, so the
ratios do not measure equal-accuracy costs for an identical nonlinear model.
CS requires $3$--$4$ corrections per step
at $5\times10^{-5}$ and $4$--$6$ at the larger steps, plus one final
residual evaluation. Mass drift remains at $10^{-13}$--$10^{-12}$.
The energy and accuracy interpretations are given in
Sections~\ref{sec:structure-tests} and~\ref{sec:production-accuracy}.

\begin{table}[t]
\centering
\small
\setlength{\tabcolsep}{4pt}
\caption{Solver comparison on the $r=4$, $p=5$ mesh (64\,002~DOFs, $\eps=0.01$, $T=2$).  SI: semi-implicit, CS: stabilized convex splitting. The reported energy is the truncated discrete energy $\mathcal E_h(T)$ of \eqref{eq:Eh}, evaluated as a common diagnostic for both schemes. The modified Lyapunov functional $E_h^n$ of \eqref{eq:Lyap} is monitored separately in Table~\ref{tab:energy_mono}. ``Terminal range'' is the range of $u_h(T)$ sampled at the vertices of a sixfold subdivided visualization mesh; ``$\max_t|u_h|$'' is the largest pointwise value over the whole trajectory, sampled at the quadrature nodes at every step (available from the semi-implicit reruns of Section~\ref{sec:defect}; the convex-splitting trajectories were not resampled). Wall times are the stepping costs of the original runs (Section~\ref{sec:environment}).}\label{tab:dt_comparison}
\begin{tabular}{lccccccc}
\toprule
$\dt$ & Solver & Steps & $\mathcal E_h(T)$ & $|\delta m|_{\max}$ & Terminal range & $\max_t|u_h|$ & Wall time \\
\midrule
\multirow{2}{*}{$5\times10^{-5}$}
  & SI & 40\,000 & 10.26 & $1.8\times10^{-12}$ & $[-1.02, 1.01]$ & $3.47$ & 22.5\,min \\
  & CS & 40\,000 & 10.21 & $1.1\times10^{-13}$ & $[-1.02, 1.01]$ & -- & 95.5\,min \\
\midrule
\multirow{2}{*}{$2\times10^{-4}$}
  & SI & 10\,000 & 9.60  & $5.0\times10^{-13}$ & $[-1.02, 1.02]$ & $3.17$ & 5.6\,min \\
  & CS & 10\,000 & 9.67  & $1.2\times10^{-13}$ & $[-1.02, 1.02]$ & -- & 28.8\,min \\
\midrule
\multirow{2}{*}{$5\times10^{-4}$}
  & SI & 4\,000  & 10.76 & $1.7\times10^{-13}$ & $[-1.08, 1.51]$ & $3.30$ & 2.2\,min \\
  & CS & 4\,000  & 10.72 & $1.0\times10^{-13}$ & $[-1.07, 1.15]$ & -- & 12.9\,min \\
\bottomrule
\end{tabular}
\end{table}

Table~\ref{tab:cs_residuals} tests the convergence of the fixed-matrix
corrections through onset. All $120$ monitored steps reach tolerance in
$2$--$6$ corrections (mean $3.5$). Contraction grows from about
$7\times10^{-6}$ to $10^{-2}$. Its scale agrees with
\eqref{eq:cs-contraction}, but it can exceed the constant-weight modal
estimate, as at step $30$. Omitting stabilization from the nonlinear
scheme while retaining the same correction matrix gives contraction
$0.9997$ and stalls the iteration.

\begin{table}[t]
\centering
\small
\caption{Convex-splitting corrections at $\dt=5\times10^{-4}$ ($r=4$, $p=5$, $\eps=0.01$, $S=600$): number of corrections (back-substitutions; the residual is evaluated once more to confirm convergence), residual norm before the first correction and at termination, observed per-correction contraction (geometric mean over the corrections after the first), the linearized estimate $3\max|u_h|^2/S$ of \eqref{eq:cs-contraction}, at representative steps through the onset window $t\in[0.024,0.039]$.}\label{tab:cs_residuals}
\begin{tabular}{ccccccc}
\toprule
step & $t$ & corrections & $r_0$ & $r_{\mathrm{end}}$ & observed contraction & $3\max|u_h|^2/S$ \\
\midrule
1 & 0.0005 & 2 & $8.94\times10^{2}$ & $1.11\times10^{-11}$ & $6.68\times10^{-6}$ & $1.25\times10^{-5}$ \\
3 & 0.0015 & 2 & $7.49\times10^{2}$ & $1.16\times10^{-11}$ & $6.80\times10^{-6}$ & $1.25\times10^{-5}$ \\
10 & 0.0050 & 2 & $3.81\times10^{2}$ & $1.71\times10^{-11}$ & $7.92\times10^{-6}$ & $1.41\times10^{-5}$ \\
30 & 0.0150 & 3 & $5.15\times10^{2}$ & $2.34\times10^{-13}$ & $1.54\times10^{-4}$ & $4.04\times10^{-5}$ \\
48 & 0.0240 & 3 & $1.12\times10^{3}$ & $6.88\times10^{-13}$ & $1.00\times10^{-4}$ & $1.53\times10^{-4}$ \\
55 & 0.0275 & 3 & $1.66\times10^{3}$ & $3.33\times10^{-12}$ & $1.48\times10^{-4}$ & $2.62\times10^{-4}$ \\
62 & 0.0310 & 3 & $1.91\times10^{3}$ & $2.15\times10^{-11}$ & $2.53\times10^{-4}$ & $4.51\times10^{-4}$ \\
70 & 0.0350 & 4 & $2.14\times10^{3}$ & $1.06\times10^{-12}$ & $1.09\times10^{-3}$ & $8.31\times10^{-4}$ \\
78 & 0.0390 & 4 & $2.66\times10^{3}$ & $1.97\times10^{-12}$ & $1.01\times10^{-3}$ & $1.51\times10^{-3}$ \\
90 & 0.0450 & 4 & $3.68\times10^{3}$ & $4.06\times10^{-11}$ & $1.97\times10^{-3}$ & $3.42\times10^{-3}$ \\
105 & 0.0525 & 5 & $4.00\times10^{3}$ & $8.09\times10^{-12}$ & $5.72\times10^{-3}$ & $7.32\times10^{-3}$ \\
120 & 0.0600 & 6 & $7.54\times10^{3}$ & $6.10\times10^{-12}$ & $1.01\times10^{-2}$ & $9.62\times10^{-3}$ \\
\bottomrule
\end{tabular}
\end{table}

\begin{figure}
\centering
\includegraphics[width=\textwidth]{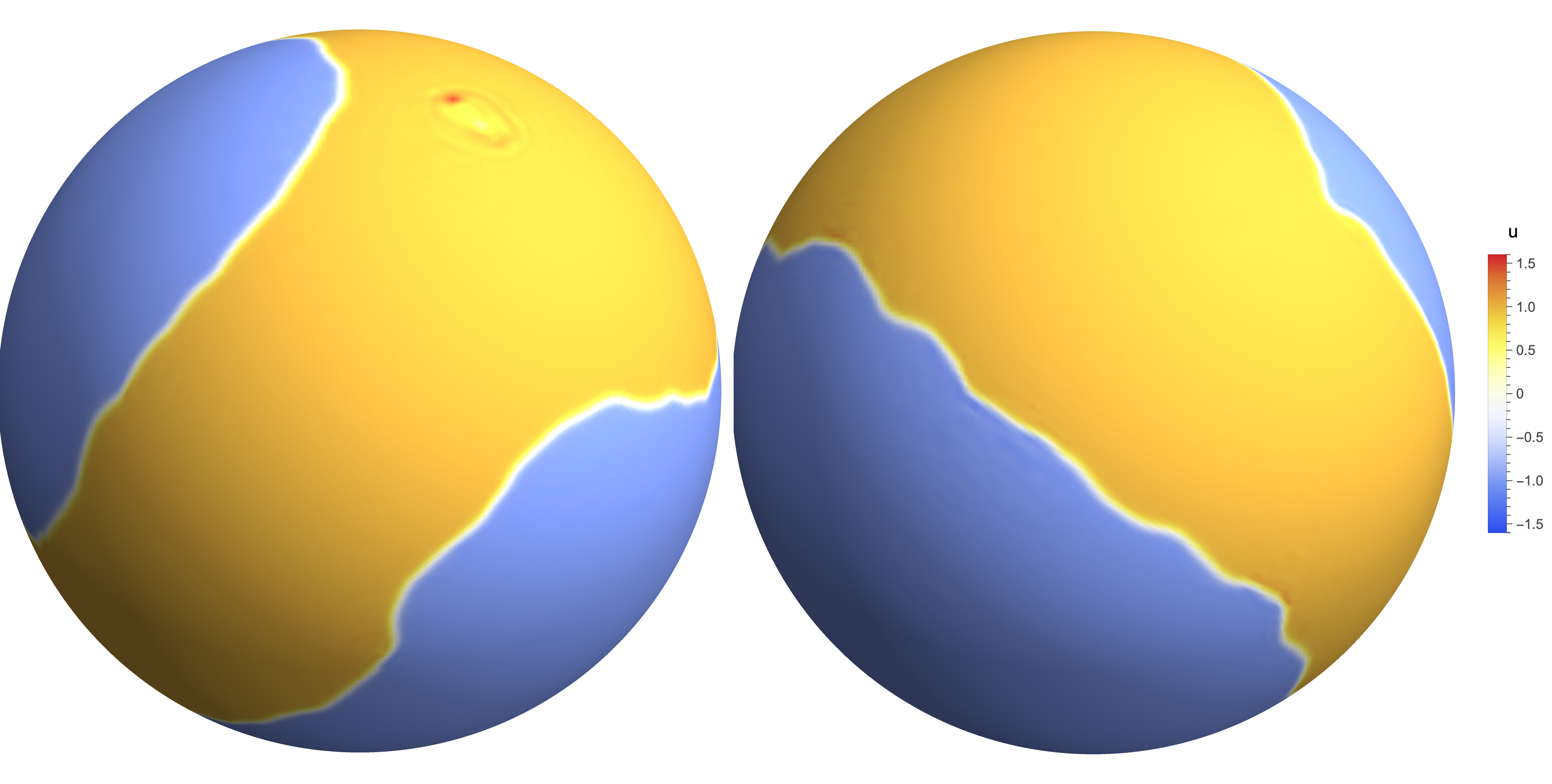}
\caption{\CapSpinodal}\label{fig:spinodal}
\end{figure}

Figure~\ref{fig:spinodal} records the terminal SI morphology at the
largest step. The accompanying table distinguishes terminal ranges from
trajectory maxima.

\subsubsection{Equal-dimension $hp$ configurations}\label{sec:hp_study}

At fixed $N_h=23\,042$ and $h/p=0.0165$, Table~\ref{tab:hp_comparison}
separates reaction assembly, solve, and initial assembly costs for SI at
$\dt=5\times10^{-4}$. Reaction assembly falls from $9.8$ to $3.5$ ms
as element count decreases. Solve cost rises from $6.1$ to $13.6$ ms as
local blocks grow. The intermediate r3$p$6 configuration has the lowest
stepping cost, while one-time assembly grows from $0.63$ to $18.76$ min.

\begin{table}[t]
\centering
\small
\setlength{\tabcolsep}{3pt}
\caption{Semi-implicit scheme at $\dt = 5\times10^{-4}$: cost breakdown across the three equal-dimension configurations ($N_h=23\,042$, $\eps=0.01$, $T=2$).  The reported energy is the truncated discrete energy $\mathcal E_h(T)$ of \eqref{eq:Eh}. $|\delta m|_{\max}$ is the maximum mass drift.  Wall time is the stepping cost, excluding assembly and the per-step diagnostics.}\label{tab:hp_comparison}
\begin{tabular}{ccccccccc}
\toprule
Config & DOFs & Assembly (min) & Factor (s) & $R_f$/step (ms) & Solve/step (ms) & Wall (min) & $\mathcal E_h(T)$ & $|\delta m|_{\max}$ \\
\midrule
r4$p$3  & 23\,042 & 0.63  & 0.44 & 9.8 & 6.1  & 1.06 & 10.97 & $2.0\times10^{-14}$ \\
r3$p$6  & 23\,042 & 2.57  & 0.66 & 3.6 & 7.2  & 0.72 & 11.34 & $4.7\times10^{-13}$ \\
r2$p$12 & 23\,042 & 18.76 & 1.11 & 3.5 & 13.6 & 1.14 & 10.79 & $3.9\times10^{-11}$ \\
\bottomrule
\end{tabular}
\end{table}

\begin{figure}
\centering
\includegraphics[width=\textwidth]{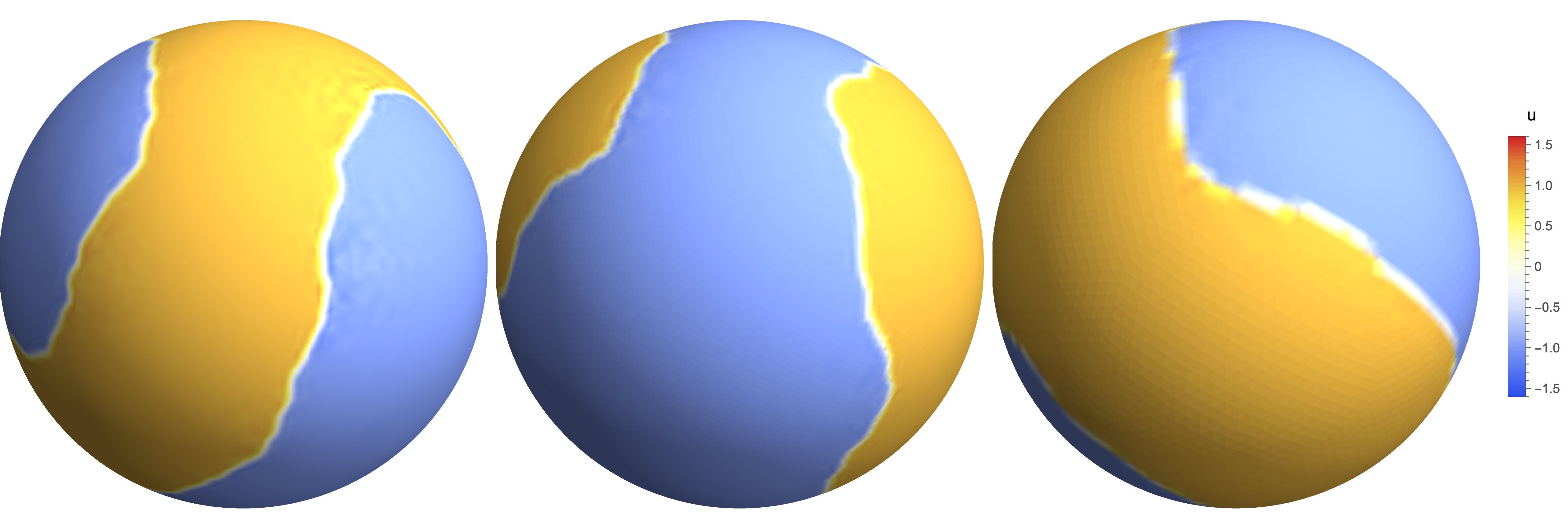}
\caption{\CapHpSnapshots}\label{fig:hp_snapshots}
\end{figure}

Mass drift rises to $3.9\times10^{-11}$ at $p=12$, reflecting solve
accuracy in the unnormalized hierarchic basis. The identity $Kc=0$ is
unchanged. Figure~\ref{fig:hp_snapshots} shows the three terminal states.
Because the initial field is projected separately into each space, these
rows compare cost at equal dimension and resolution scale. A comparison
at equal accuracy would require errors against a common converged state.

\FloatBarrier
\subsection{Geometry and topology benchmarks}\label{sec:geometry-tests}

Controlled torus bands compare diffuse energies with analytical interface
lengths. Coarsening runs on the torus and implicit higher-genus surfaces
test geometry-dependent patterns, winding classes, and structure
preservation. A viscous torus model extends the benchmark set.

\subsubsection{Torus interface energies and winding}\label{sec:torus}\label{sec:torus-bands}

Controlled torus bands test the relation between diffuse energy and
metric-dependent interface length. The chart is
\[
X(\theta,\phi)=\bigl((R+r\cos\phi)\cos\theta,
(R+r\cos\phi)\sin\theta,r\sin\phi\bigr),
\qquad (\theta,\phi)\in[0,2\pi)^2.
\]
For $R=1$, $r=1/2$, the quartic line tension and two equal-area
boundary lengths are
\begin{equation}
\sigma_0=\int_{-1}^{1}\sqrt{2F(s)}\,\dd s=\frac{2\sqrt2}{3}.
\label{eq:line-tension}
\end{equation}
\begin{equation}
L_P=4\pi r.
\label{eq:LP}
\end{equation}
\begin{equation}
\frac{R\delta-r\sin\delta}{\pi R}=\frac12,
\qquad
L_T=4\pi\bigl(R-r\cos\delta\bigr).
\label{eq:LT}
\end{equation}
where $P$ and $T$ denote poloidal and toroidal boundary loops.
The half-width is $\delta=2.0209799$, and $L_T/L_P=2.4351309$.

Profiles are $L^2$ projections of $\tanh(d/(\sqrt2\eps))$, with signed
distance positive inside the band. For poloidal boundaries centred at
$\theta_c$, $d=(\pi/2-|\theta-\theta_c|)(R+r\cos\phi)$ along
$\theta$-lines, using the periodic angular difference. For toroidal
boundaries centred at $\phi=\pi$, $d=r(\delta-|\phi-\pi|)$ along
meridians. The mesh is a structured $32\times16$ periodic grid split
into $1\,024$ triangles, with $p=5$, $N_h=12\,800$, and $\eps=0.025$.

Orientation is measured by
\begin{equation}
A_{\mathrm{tube}}:=|\hat u_{1,0}|^2,\qquad
A_{\mathrm{hole}}:=|\hat u_{0,1}|^2,\qquad
\hat u_{mn}:=\frac{1}{|\Gamma|}\int_\Gamma u_h\,e^{-i(m\theta+n\phi)}\,\dd S,
\label{eq:torus-fourier}
\end{equation}
and $\chi=(A_{\mathrm{tube}}-A_{\mathrm{hole}})/
(A_{\mathrm{tube}}+A_{\mathrm{hole}})$. For winding, values on a periodic
$128\times64$ grid are thresholded at the mean. Breadth-first search
with four-neighbour adjacency and integer lifts detects whether the
largest component of each phase wraps the $\theta$ cycle, the $\phi$
cycle, both, or neither.

\begin{table}[pos=!htbp]
\centering
\caption{Controlled equal-area torus bands. The relative error compares
$\mathcal E_h(u_h^0)/\sigma_0$ with the corresponding perimeter. The
orientation index $\chi$ distinguishes the two band directions.}
\label{tab:torus-bands}
\begin{tabular}{lccccc}
\toprule
boundary class & $L$ & $\mathcal E_h(u_h^0)/\sigma_0$ & rel. error &
$\chi(0)$ & $|\delta m|_{\max}$ \\
\midrule
poloidal & 6.28319 & 6.30045 & $2.75\times10^{-3}$ & 1 & $6.7\times10^{-14}$ \\
toroidal & 15.30038 & 15.44889 & $9.71\times10^{-3}$ & $-1$ & $5.5\times10^{-14}$ \\
\bottomrule
\end{tabular}
\end{table}

The diffuse energies reproduce the perimeters within $1\%$. Their ratio
$2.45203$ differs from $L_T/L_P$ by $0.69\%$. Short SI evolutions with
$S=2$, $\dt=2\times10^{-4}$, and $T=0.05$ retain the winding classes
and $\chi=\pm1$. Energy decreases at the eleven recordings (every $25$
steps), apart from one roundoff increment of $1.6\times10^{-14}$.

\subsubsection{Torus coarsening and orientation changes}

\begin{figure}[pos=ht]
\centering
\includegraphics[width=\textwidth]{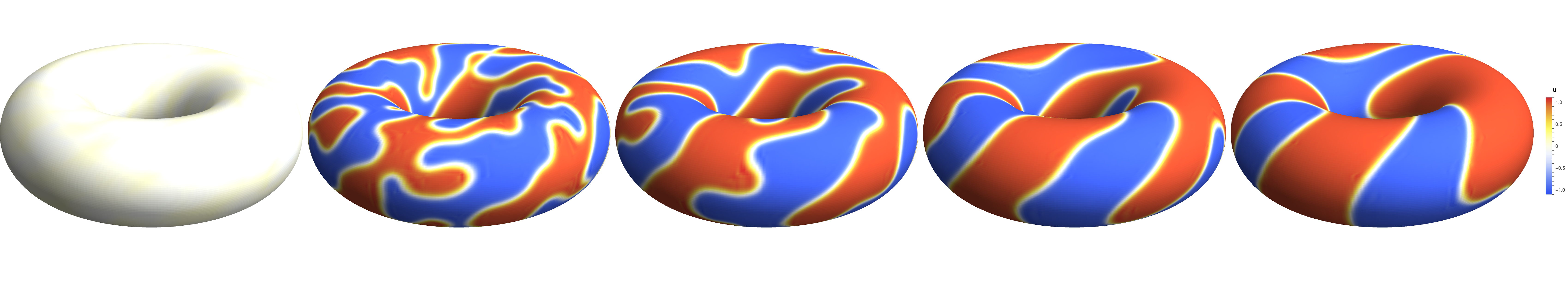}
\caption{\CapTorusPlain}\label{fig:torus_plain}
\end{figure}

Figure~\ref{fig:torus_plain} gives a random-seed SI trajectory on
$R=1$, $r=1/2$, with $\eps=0.025$, $S=12$, $\dt=5\times10^{-4}$,
$p=5$, $N_h=12\,800$, and $T=1.5$. The initial perturbation has zero
mean, and the final state contains bands along the torus cycles.

To make the orientation test reproducible, a second run uses $R=1$,
$r=0.35$, and
\begin{equation}
u_h^0=P_h\left[0.05\bigl(\cos\theta+\cos\phi
+0.2\cos(2\theta-\phi)\bigr)\right].
\label{eq:torus-seed}
\end{equation}
Its mean is $\bar u=0.00875$, giving phase fraction $0.504375$.
The two-class sharp-interface lengths are $L_{\mathrm{tube}}=4.39823$
and $L_{\mathrm{hole}}=14.04688$. The SI settings are $\eps=0.025$,
$S=2$, $\dt=5\times10^{-4}$, $T=2$, $p=5$, and a $40\times20$
grid ($1\,600$ triangles, $N_h=20\,000$).

\begin{table}[pos=!htbp]
\centering
\caption{Deterministic torus benchmark \eqref{eq:torus-seed}. The
orientation index is
$\chi=(A_{\mathrm{tube}}-A_{\mathrm{hole}})/
(A_{\mathrm{tube}}+A_{\mathrm{hole}})$.}
\label{tab:torus-deterministic}
\begin{tabular}{cccc}
\toprule
$t$ & $A_{\mathrm{tube}}$ & $A_{\mathrm{hole}}$ & $\chi$ \\
\midrule
0.05 & $9.86\times10^{-4}$ & $2.09\times10^{-2}$ & $-0.910$ \\
0.20 & $7.31\times10^{-4}$ & $2.46\times10^{-2}$ & $-0.942$ \\
0.50 & $4.50\times10^{-3}$ & $2.30\times10^{-2}$ & $-0.672$ \\
2.00 & $1.34\times10^{-1}$ & $1.73\times10^{-3}$ & $0.975$ \\
\bottomrule
\end{tabular}
\end{table}

Fourier content changes from hole-dominant to tube-dominant, while at
$T=2$ the largest connected component of each phase is classified as
wrapping both periodic directions. Energy falls from $137.4744$ to
$10.9695$ at recordings every
ten steps. The largest recorded increment is $-4.8\times10^{-3}$,
the ranges at the four times in Table~\ref{tab:torus-deterministic}
lie within $[-1.01,1.01]$, and the recorded maximum of $|u_h|$ over
the run is $1.30$. Mass drift is $1.7\times10^{-13}$. The final wrapping
configuration lies outside the two coordinate-aligned band classes
used in the perimeter comparison. Their energy ordering does not establish
equilibrium of this final configuration. These observations concern the stated
finite $\eps$, mesh, and step. An $\eps$ limit would require independent
refinement of $\eps$, $h$, and $\dt$.

\subsubsection{Viscous surface-model extension}\label{sec:voigt-test}

The R{\"a}tz--Voigt benchmark tests a viscous extension on the
$R=1$, $r=1/2$ torus:
\begin{align}
\partial_t u &= \nu\Delta_\Gamma \mu,
\label{eq:voigt-1}
\\
\mu &= -\gamma\Delta_\Gamma u
      +\gamma^{-1}G'(u)
      +\alpha\gamma\,\partial_t u,
\label{eq:voigt-2}
\end{align}
\[
G(u)=18u^2(1-u)^2,\qquad \gamma=0.05,\quad\alpha=\nu=1,
\qquad
\mathcal E_V=\frac\gamma2\int_\Gamma|\nabla_\Gamma u|^2\,\dd S
+\gamma^{-1}\int_\Gamma G(u)\,\dd S.
\]
Viscosity adds the continuous dissipation
$-\alpha\gamma\norm{\partial_tu}_{L^2(\Gamma)}^2$.
The discrete energy analysis of Theorem~\ref{thm:lyap} applies to the
non-viscous semi-implicit scheme.

The implementation modifies SAV-BDF2 \citep{ShenXuYang2018}. With
$r=\sqrt{\int_\Gamma G(u)\,\dd S+C}$, $C=1$, and
$\bar u=2u^n-u^{n-1}$, the reaction is replaced by
\[
\gamma^{-1}\frac{r^{n+1}G'(\bar u)}
{\sqrt{\int_\Gamma G(\bar u)\,\dd S+C}}.
\]
The auxiliary variable follows its BDF2 discrete chain rule, viscosity
uses the implicit BDF2 difference, and $S(u^{n+1}-\bar u)$ is added
with $S=360=\max_{[0,1]}|G''|/(2\gamma)$. The quartic has unbounded
$|G''|$ on $\mathbb R$, and the trajectory leaves $[0,1]$ slightly. After
each step, $r^{n+1}$ is reset to its defining integral. This addresses
auxiliary-variable consistency, for which \citet{JiangZhangZhao2022}
propose constrained relaxation. The unconditional reset used here is
assessed empirically through the observed decrease of $\mathcal E_V$.
An SAV energy inequality for this reset would require a separate proof.

\begin{figure}[pos=ht]
\centering
\includegraphics[width=\textwidth]{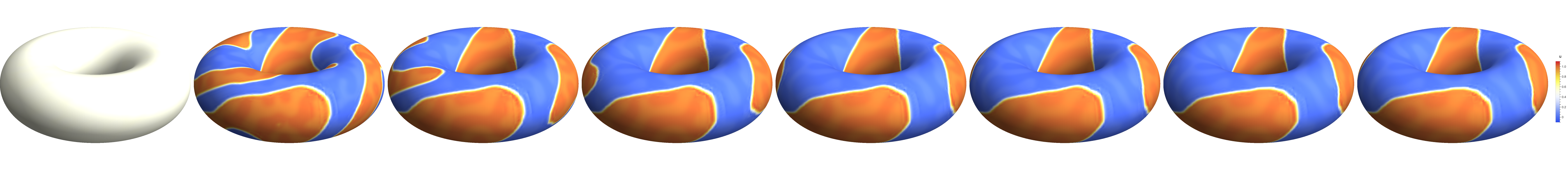}
\caption{\CapTorusVoigt}\label{fig:torus_voigt}
\end{figure}

Two backward-Euler half-steps start the run. A sparse multifrontal solve
handles the bordered $(2N_h+1)$ system. The mesh has $4\,096$ triangles,
$p=3$, and $N_h=18\,432$. The time step is $\dt=10^{-3}$ and $T=1$. Initial values
perturb the mean $0.5$ by amplitude $0.01$.
Table~\ref{tab:voigt_torus} and Figure~\ref{fig:torus_voigt} give the
mass, physical-energy change, pointwise range, and coarsening sequence.

\begin{table}[pos=ht]
\centering
\caption{Viscous R{\"a}tz--Voigt torus coarsening at $T=1$ with $p=3$ and $\dt=10^{-3}$.  The reported energy change uses $\mathcal E_V$.}\label{tab:voigt_torus}
\begin{tabular}{ccccc}
\toprule
DOFs & $|\delta m|_{\max}$ & $\Delta\mathcal E_V$ & $u_h$ range & final time \\
\midrule
18\,432 & $7.5\times10^{-14}$ & $-4.23\times10^2$ & $[-0.092,\,1.092]$ & $1.0$ \\
\bottomrule
\end{tabular}
\end{table}

\FloatBarrier
\subsubsection{Implicit higher-genus geometries}\label{sec:highergenus}

These tests assess the method on closed surfaces represented by implicit
level sets. The genus-2 and genus-3 benchmarks of \citet{hrz}
are the level sets $\Phi_2=13/10$ and $\Phi_3=13/10$, where
\begin{equation}
\Phi_2(x,y,z)
= \frac{3}{32}\bigl(x^2+2y^2\bigr)^2
+ \frac{1}{2\bigl((x-\tfrac34)^2+y^2\bigr)+\tfrac12}
+ \frac{1}{2\bigl((x+\tfrac34)^2+y^2\bigr)+\tfrac12}
+ 2z^2
\label{eq:phi-genus2}
\end{equation}
\begin{equation}
\Phi_3(x,y,z)
= \frac{1}{24}\bigl(x^2+y^2\bigr)^2
+ \frac{2}{4x^2+4(y-2)y+5}
+ \frac{2}{4x(x-\sqrt3)+4y(y+1)+5}
+ \frac{2}{4x(x+\sqrt3)+4y(y+1)+5}
+ z^2
\label{eq:phi-genus3}
\end{equation}
A triangular seed mesh is extracted and its vertices projected by
$\mathbf x\leftarrow\mathbf x-(\Phi-c)\nabla\Phi/|\nabla\Phi|^2$
to residual $10^{-12}$. On each element, $\mathbf x_K=P\circ Y_K$
composes the affine seed map with closest-point projection. At quadrature
points, Newton iteration solves
\[
\mathbf x+\lambda\nabla\Phi(\mathbf x)=\mathbf y,\qquad
\Phi(\mathbf x)=c,\qquad
J=\begin{bmatrix}I+\lambda\nabla^2\Phi&\nabla\Phi\\
\nabla\Phi^T&0\end{bmatrix},
\]
starting at $(\mathbf y,0)$ with tolerance $10^{-12}$. The same Jacobian
gives the tangents $DP(Y_K)\partial_\alpha Y_K$. Thus quadrature uses
the projected-surface metric, with the degrees specified in
Section~\ref{sec:spatial}. Manifoldness, projection residuals, and Euler
characteristics $736-2214+1476=-2$ and $1384-4164+2776=-4$ verify the meshes.

\begin{figure}[pos=!htbp]
\centering
\includegraphics[width=\textwidth]{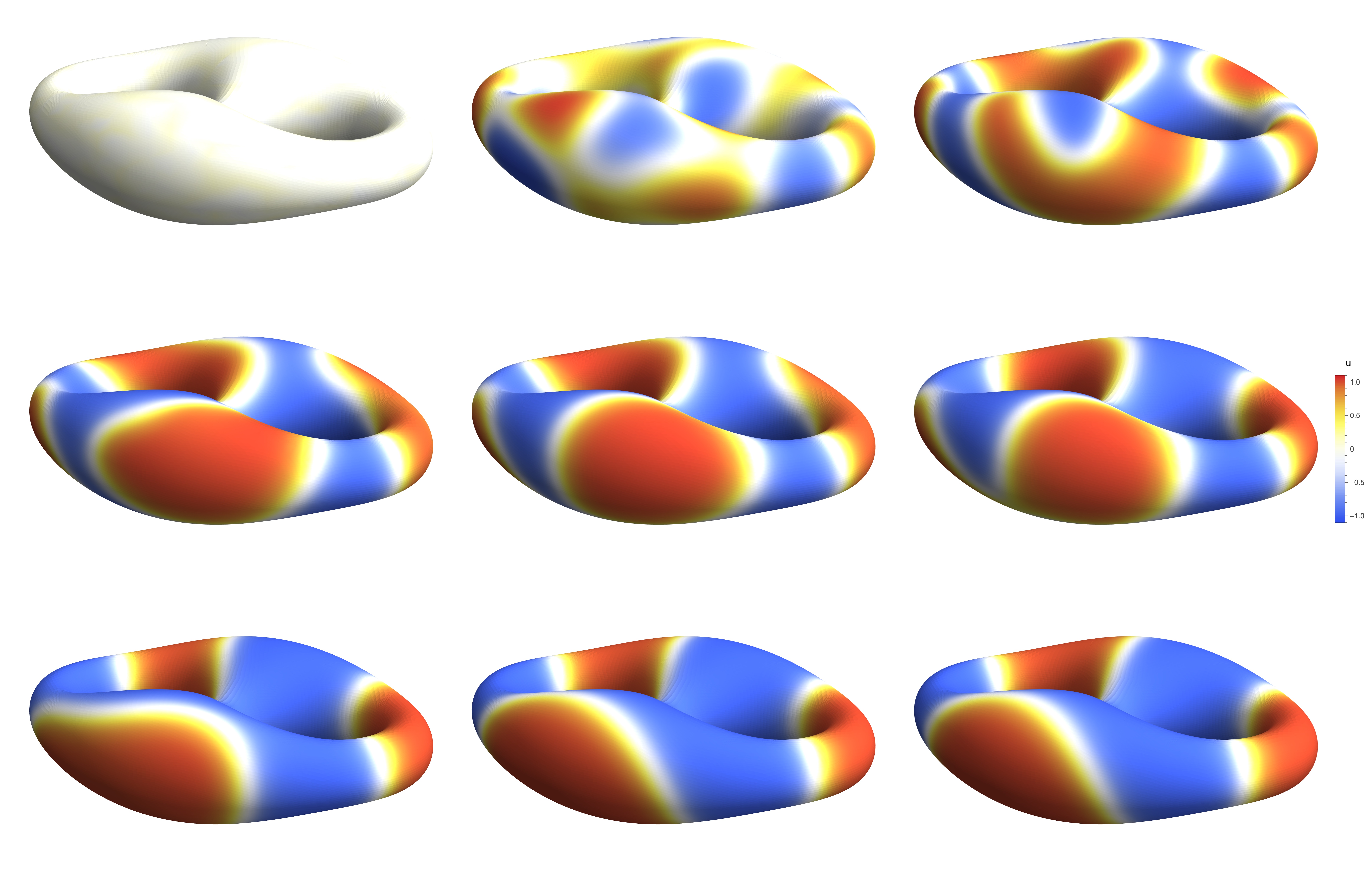}
\caption{\CapGenusTwoGrid}\label{fig:genus2_grid}
\end{figure}

\begin{figure}[pos=!htbp]
\centering
\includegraphics[width=\textwidth]{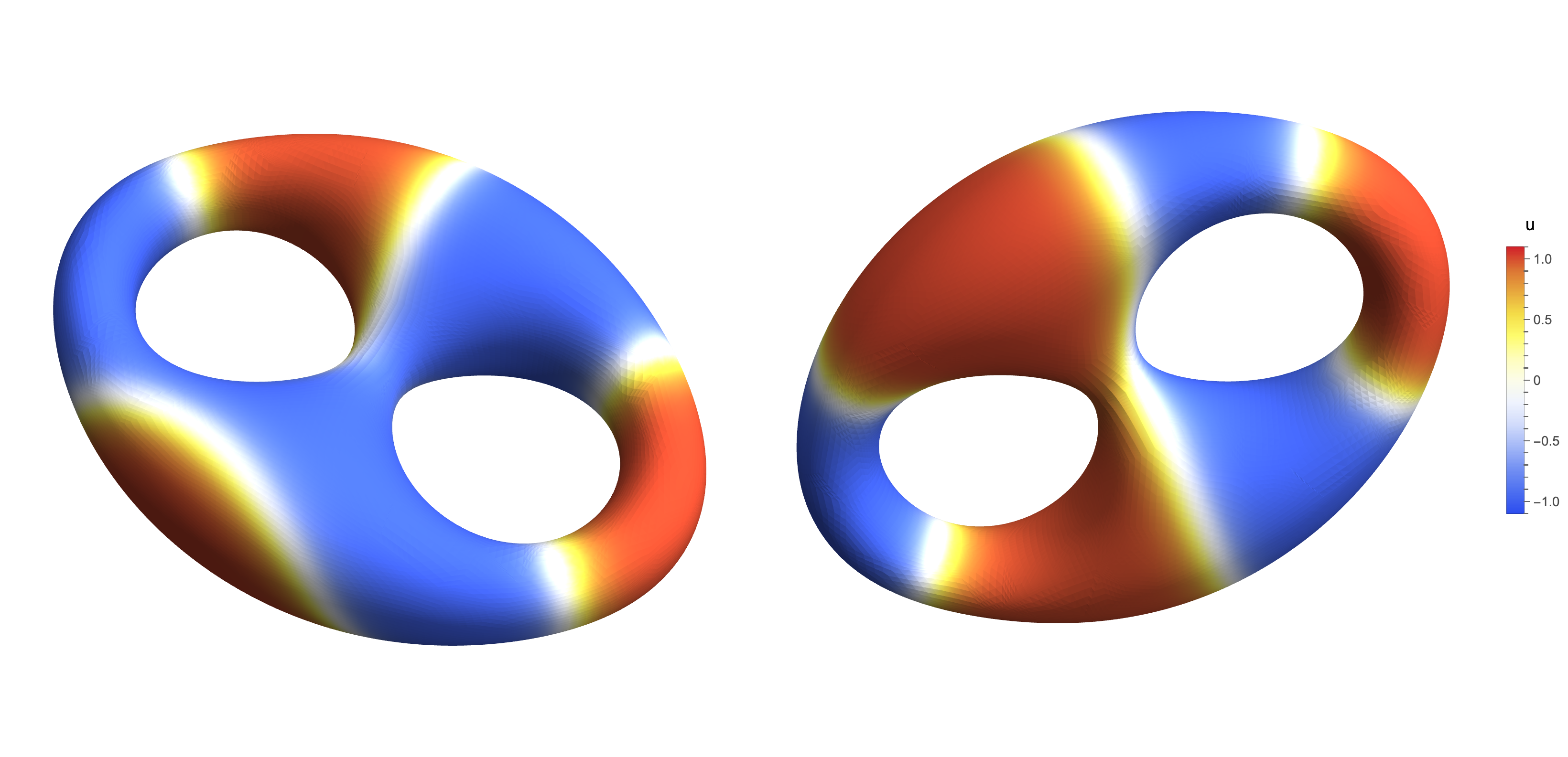}
\caption{\CapGenusTwoHero}\label{fig:genus2_hero}
\end{figure}

\begin{figure}[pos=!htbp]
\centering
\includegraphics[width=\textwidth]{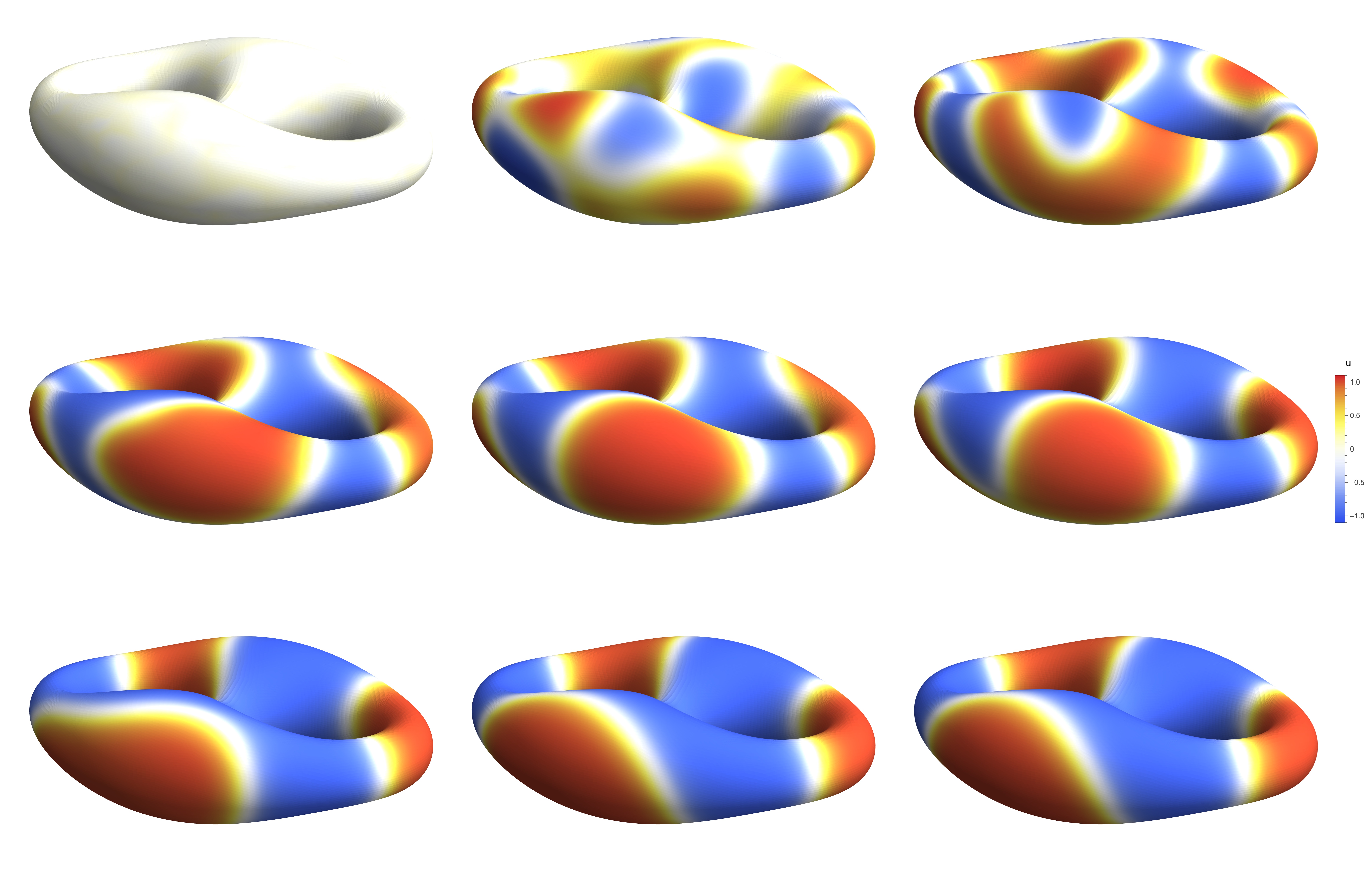}
\caption{\CapGenusThreeGrid}\label{fig:genus3_grid}
\end{figure}

\begin{figure}[pos=!htbp]
\centering
\includegraphics[width=\textwidth]{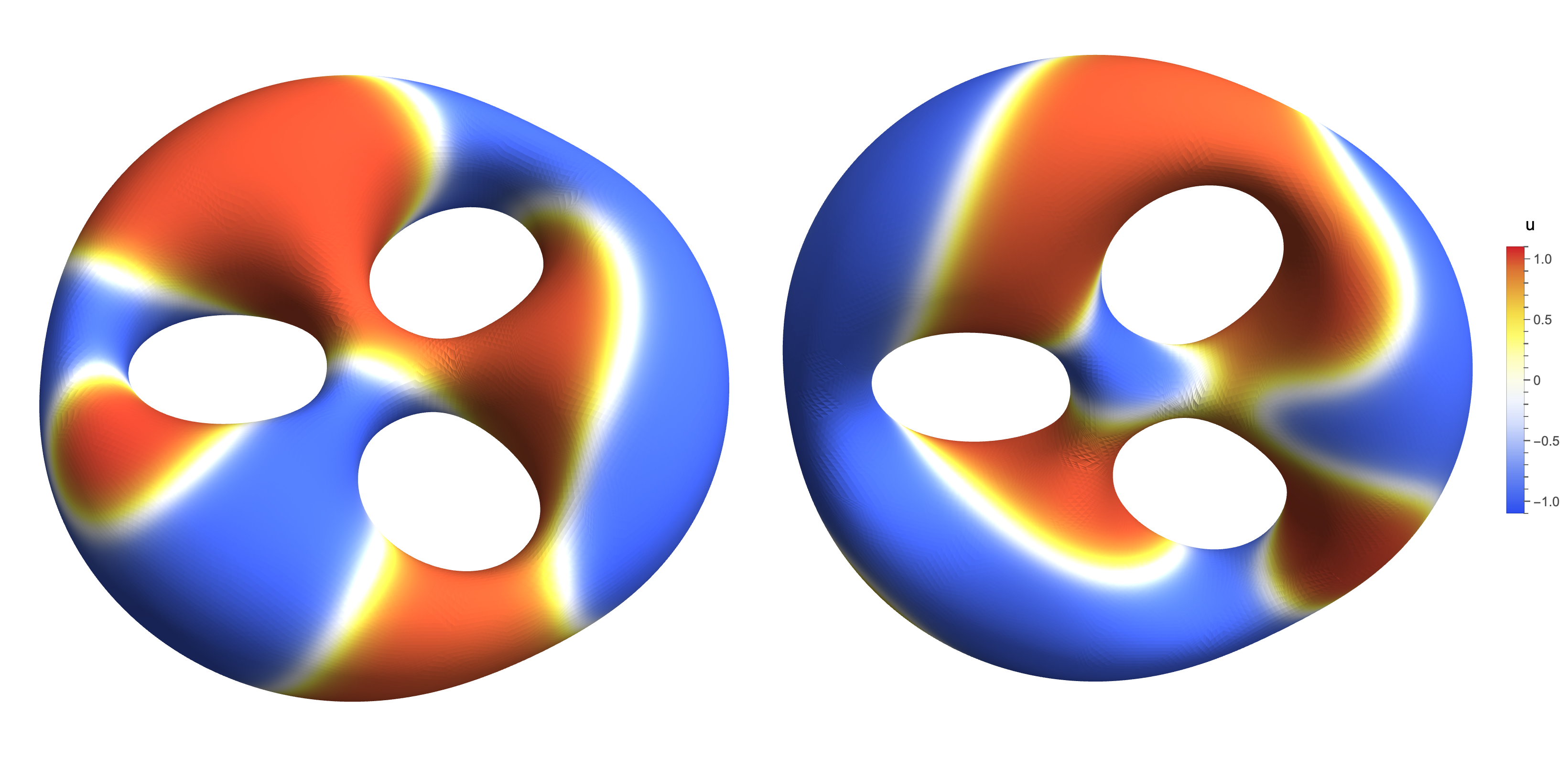}
\caption{\CapGenusThreeHero}\label{fig:genus3_hero}
\end{figure}

Both runs use SI with $p=3$, $h=0.3$, $\eps=0.1$, $S=12$,
$\dt=2\times10^{-3}$, $T=1$, and a zero-mean random vertex perturbation
of amplitude $0.1$. Table~\ref{tab:genus23} reports energy changes and
pointwise diagnostics. Figures~\ref{fig:genus2_grid}--\ref{fig:genus3_hero}
show the trajectories and final states.

\begin{table}[pos=ht]
\centering
\caption{Cahn--Hilliard coarsening on the genus-2 and genus-3 implicit
surfaces ($p=3$, $h=0.3$, $\eps=0.1$, $S=12$, $\dt=2\times10^{-3}$,
$T=1$, semi-implicit BDF2, 500 steps).  $|\delta m|_{\max}$ is the
maximum mass drift over the run. $\norm{u_h}_{L^\infty}^{\max}$ is the peak
pointwise value of $u_h$ evaluated at the quadrature points, and ``trunc.'' is
the maximum truncation-active area fraction.  Both $\mathcal E_h$ and $E_h^n$
are monotone at every step in both runs.}\label{tab:genus23}
\begin{tabular}{lccccccccc}
\toprule
Surface & $V/E/F$ & DOFs & Area & $\mathcal E_h(0)$ & $\mathcal E_h(T)$ &
Drop & $|\delta m|_{\max}$ & $\norm{u_h}_{L^\infty}^{\max}$ & trunc.\ (\%) \\
\midrule
Genus 2 & 736/2214/1476  & 6\,640  & 15.976 & 40.575  & 11.877 & 70.7\,\% & $5.5\times10^{-14}$ & 1.29 & 15.3 \\
Genus 3 & 1384/4164/2776 & 12\,488 & 39.958 & 100.932 & 24.677 & 75.6\,\% & $2.6\times10^{-14}$ & 1.20 & 13.4 \\
\bottomrule
\end{tabular}
\end{table}

Both energies decrease at all $500$ steps, with mass drift at
$10^{-14}$. Truncation is active: pointwise peaks $1.29$ and $1.20$
occur during onset, with maximal active-area fractions $15.3\%$ and
$13.4\%$, and terminal maxima $1.01$ and $1.02$.
At $\dt=10^{-3}$ the energy drops become $64.3\%$ and $75.4\%$
(from $70.7\%$ and $75.6\%$), peaks fall to $1.04$ and $1.05$,
and active fractions to $11.7\%$ and $12.7\%$. The genus-2 trajectory
is appreciably time-step sensitive at the coarser step.

The unstable band has $0<\lambda<100$. The discrete counts of $126$
and $317$ nonconstant modes agree within $1\%$ with the area-based
Weyl estimates $A\lambda_c/(4\pi)=127$ and $318$. Their ratio $2.52$
is close to the area ratio $2.50$. The surfaces also differ in metric,
neck geometry, and initial spectral content. Isolating genus dependence
in the energy drops requires controlling these features.

\subsubsection{Interpretation of geometry and topology}\label{sec:representation}

Analytic charts supply the metric on the sphere and torus. Projection
derivatives supply it on the implicit surfaces. Constants belong to each
space and $Kc=0$ supplies the same conservation mechanism. Different
dimensions and geometries preclude a direct timing comparison of the
two representations.

The metric determines interface lengths and modal spectra. Topology
constrains the available winding classes: the complete phase boundary is
null-homologous, so its component homology classes cancel. On a torus,
a noncontractible simple loop is non-separating and a wrapping band has
paired boundary loops. For genus at least two, a separating
noncontractible loop can itself be null-homologous and bound a phase
domain alone. The displayed higher-genus states are bands with paired
loops. Their diffuse energy admits the line-tension interpretation
\eqref{eq:line-tension}. The present model contains no explicit
composition--curvature coupling, although geometry can bias domain
organization through the surface metric. Explicit curvature-dependent
chemical preferences can be incorporated on prescribed surfaces, whereas
budding and reciprocal shape--composition feedback additionally require
surface deformation
\citep{BaumgartHessWebb2003,CallanJonesSorreBassereau2011,JulicherLipowsky1993}.
\FloatBarrier

\section{Conclusion}\label{sec:conclusion}

We studied a mixed $hp$-finite element discretization for the Cahn--Hilliard equation on curved closed surfaces. The computations use the unit sphere, a torus, and implicit level-set surfaces of genus two and three. The time integrator is the stabilized semi-implicit BDF2 scheme with vector reaction extrapolation, which is shared with earlier work on the sphere \citep{ZhangWangRong2011}. The contribution lies in the $hp$ surface discretization on curved elements, in the discrete structure statements that hold algebraically for it, in the diagnostics that test their one non-elementary hypothesis, and in the geometry benchmarks. A stabilized convex-splitting BDF2 variant, solved by fixed-matrix corrections with the semi-implicit factorization, serves as the comparison.

On the theoretical side, we proved exact discrete mass conservation (Theorem~\ref{thm:mass}), energy decay of the truncated discrete energy $\mathcal E_h$ for the stabilized start-up steps (Lemma~\ref{lem:startup}), an exact per-step energy balance for the BDF2 phase (Corollary~\ref{cor:balance}), and a conditional monotone law for a modified discrete Lyapunov functional $E_h^n$ (Theorem~\ref{thm:lyap}). The balance isolates the extrapolation defect as its sole sign-indefinite term. The diagnostic of Section~\ref{sec:defect} evaluates it along computed trajectories and verifies the balance to roundoff.

A manufactured solution with an analytic forcing, chosen on a linearly stable state, verifies second-order convergence in time and the expected $h$- and $p$-convergence of the mixed discretization. The single-mode sphere test converges to the analytical rate \eqref{eq:growth-rate}. Across the production runs at $\eps=0.01$, semi-implicit mass drift is $10^{-14}$--$10^{-11}$, and the convex-splitting runs satisfy the invariant of Theorem~\ref{thm:cs-mass} to the same scale. The modified energy $E_h^n$ decreases at every step for the two larger time steps. At the smallest time step, the original run shows late increases of at most $1.1\times10^{-7}$ (Table~\ref{tab:energy_mono}). The defect diagnostic identifies the onset interval where the sufficient bound in Assumption~\ref{ass:ext} is exceeded for $S=600$. The exact balance accounts for the decrease of $E_h^n$ during that interval. The convex-splitting corrections contract at a rate consistent in scale with the linearized estimate $3\bar u^2/S$ of \eqref{eq:cs-contraction}. At fixed mesh and time step, its elapsed-time cost is $4.2$--$5.8$ times that of the semi-implicit implementation. The two schemes discretize different potentials outside $[-1,1]$.

The production time steps exceed the fastest linear time scale by one to two orders of magnitude. A single-mode test at production parameters shows that stabilization flattens the dispersion relation and numerical damping governs onset. A post-onset refinement study measures the time-step error during coarsening (Section~\ref{sec:production-accuracy}). The production runs illustrate relaxation with time-step-dependent onset kinetics and terminal states.

At fixed global dimension, the three $hp$ configurations have stepping
times between $0.72$ and $1.14$ minutes. Reaction assembly decreases with
$p$, while solve time, one-time assembly cost, and mass-matrix conditioning
increase. The r3$p$6 configuration has the smallest stepping time in
Table~\ref{tab:hp_comparison}.

On the torus, controlled diffuse energies reproduce two analytical
interface perimeters within $1\%$. The deterministic benchmark evolves from
hole-dominant to tube-dominant Fourier content. At $T=2$, the largest
connected component of each phase is classified as wrapping both periodic
directions. The viscous R{\"a}tz--Voigt run, computed with a stabilized SAV-BDF2 integrator, conserves mass to
roundoff and decreases its reported energy. On the genus-two and genus-three
surfaces, mass drift remains at $10^{-14}$ and both monitored energies are
monotone at every step. Their unstable-mode counts agree with area-based
Weyl estimates within $1\%$.

These results separate two effects of geometry. The surface metric sets the
modal spectrum and interface lengths. Topology sets the available winding
classes. The calculations establish these metric and topological effects for
the passive model. Geometry can bias domain organization through the
surface metric. Explicit curvature-dependent chemical preferences can be
incorporated on prescribed surfaces through composition--curvature
coupling. Budding and reciprocal shape--composition feedback additionally
require surface deformation. A genus-dependent coarsening law requires
controlled geometry and initial spectral content.

A direction for future work is to extend the present $hp$-surface finite
element method to the sharp-surface multiphysics models discussed by
\citet{Xu2026OVPDDM}, coupling phase separation to membrane deformation,
hydrodynamics, and active fields.

\section*{CRediT authorship contribution statement}
\textbf{Harri Hakula:} Conceptualization, Methodology, Software, Validation, Formal analysis, Investigation, Writing -- original draft, Writing -- review \& editing, Visualization. 
\textbf{Antti Rasila:} Conceptualization, Methodology, Formal analysis, Writing -- review \& editing. 
\textbf{Xinpeng Xu:} Conceptualization, Methodology, Formal analysis, Validation, Writing -- review \& editing.

\section*{Declaration of competing interest}
The authors declare that they have no known competing financial interests or personal relationships that could have appeared to influence the work reported in this paper.

\section*{Data and code availability}
The finite element framework used in this study is under active development.
The framework and the numerical data supporting the reported results are
available from the corresponding author on request.

\section*{Funding}
A.R. is supported by the Natural Science Foundation of Guangdong Province
(No.~2024A1515010467) and the Li Ka Shing Foundation GTIIT--STU Joint
Research Grant (No.~2024LKSFG06).
X.X. is supported by the National Natural Science Foundation of China
(NSFC, Nos. 12131010 and 12374209).
H.H. is supported by the Research Council of Finland 
(Flagship of Advanced Mathematics for Sensing Imaging and Modelling grant 
359181).

\section*{Declaration of generative AI and AI-assisted technologies in the writing process}
During the preparation of this work the authors used OpenAI Codex to
assist with manuscript language revision and consistency checks.
AI tools have been used in development of the underlying FEM framework,
in particular regression testing.
The authors have reviewed and edited the content and take full responsibility for the content of the publication.

\bibliographystyle{cas-model2-names}
\bibliography{references}

@misc{Xu2026OVPDDM,
  author        = {Xu, Xinpeng},
  title         = {Onsager-variational formulation of diffuse-domain methods for computational modeling of microscale fluid-structure interactions},
  year          = {2026},
  eprint        = {2605.13196},
  archivePrefix = {arXiv},
  primaryClass  = {physics.bio-ph},
  doi           = {10.48550/arXiv.2605.13196}
}

@article{CahnHilliard1958,
  author = {Cahn, J. W. and Hilliard, J. E.},
  title = {Free energy of a nonuniform system. I. Interfacial free energy},
  journal = {Journal of Chemical Physics},
  volume = {28},
  number = {2},
  pages = {258--267},
  year = {1958},
  doi = {10.1063/1.1744102}
}

@article{CondetteMelcherSuli2011,
  author = {Condette, N. and Melcher, C. and S{\"u}li, E.},
  title = {Spectral approximation of pattern-forming nonlinear evolution equations with double-well potentials of quadratic growth},
  journal = {Mathematics of Computation},
  volume = {80},
  number = {273},
  pages = {205--223},
  year = {2011},
  doi = {10.1090/S0025-5718-10-02365-3}
}

@article{LiKimWang2017,
  author  = {Li, Yibao and Kim, Junseok and Wang, Nan},
  title   = {An Unconditionally Energy-Stable Second-Order Time-Accurate Scheme
             for the {Cahn--Hilliard} Equation on Surfaces},
  journal = {Communications in Nonlinear Science and Numerical Simulation},
  volume  = {53},
  pages   = {213--227},
  year    = {2017},
  doi     = {10.1016/j.cnsns.2017.05.006}
}

@article{MarenduzzoOrlandini2013,
  author  = {Marenduzzo, Davide and Orlandini, Enzo},
  title   = {Phase Separation Dynamics on Curved Surfaces},
  journal = {Soft Matter},
  volume  = {9},
  pages   = {1178--1187},
  year    = {2013},
  doi     = {10.1039/C2SM27081A}
}

@article{ElliottRanner2015,
  author = {Elliott, C. M. and Ranner, T.},
  title = {Evolving surface finite element method for the {Cahn--Hilliard} equation},
  journal = {Numerische Mathematik},
  volume = {129},
  number = {3},
  pages = {483--534},
  year = {2015},
  doi = {10.1007/s00211-014-0644-y}
}

@inproceedings{Eyre1998,
  author = {Eyre, D. J.},
  title = {Unconditionally gradient stable time marching the {Cahn--Hilliard} equation},
  booktitle = {Computational and Mathematical Models of Microstructural Evolution},
  series = {MRS Proceedings},
  volume = {529},
  pages = {39--46},
  year = {1998},
  doi = {10.1557/PROC-529-39}
}

@incollection{Shen2011,
  author    = {Shen, Jie},
  title     = {Modeling and Numerical Approximation of Two-Phase Incompressible Flows by a Phase-Field Approach},
  booktitle = {Multiscale Modeling and Analysis for Materials Simulation},
  editor    = {Bao, Weizhu and Du, Qiang},
  series    = {Lecture Notes Series, Institute for Mathematical Sciences, National University of Singapore},
  volume    = {22},
  pages     = {147--195},
  publisher = {World Scientific},
  year      = {2011},
  doi       = {10.1142/9789814360906_0003}
}

@article{WangYu2018,
  author  = {Wang, Lin and Yu, Haijun},
  title   = {On Efficient Second Order Stabilized Semi-Implicit Schemes for the {Cahn--Hilliard} Phase-Field Equation},
  journal = {Journal of Scientific Computing},
  volume  = {77},
  number  = {2},
  pages   = {1185--1209},
  year    = {2018},
  doi     = {10.1007/s10915-018-0746-2}
}

@article{hrz,
  author  = {Hakula, Harri and Rasila, Antti and Zheng, Yufan},
  title = {Efficient Numerical Conformal Mappings on Multiply Connected {R}iemann Surfaces},
  journal = {Electronic Transactions on Numerical Analysis},
  volume  = {65},
  pages   = {533--552},
  year    = {2026}
}

@article{LiaoJiWangZhang2022,
  author = {Liao, H.-L. and Ji, B. and Wang, L. and Zhang, Z.},
  title = {Mesh-robustness of an energy stable {BDF2} scheme with variable steps for the {Cahn--Hilliard} model},
  journal = {Journal of Scientific Computing},
  volume = {92},
  number = {2},
  pages = {52},
  year = {2022},
  doi = {10.1007/s10915-022-01861-4}
}

@article{YanChenWangWise2018,
  author  = {Yan, Y. and Chen, W. and Wang, C. and Wise, S. M.},
  title   = {A second-order energy stable {BDF} numerical scheme for the {Cahn--Hilliard} equation},
  journal = {Communications in Computational Physics},
  volume  = {23},
  number  = {2},
  pages   = {572--602},
  year    = {2018},
  doi     = {10.4208/cicp.OA-2016-0197}
}

@article{RatzVoigt2006,
  author = {R{\"a}tz, A. and Voigt, A.},
  title = {{PDE}s on surfaces: a diffuse interface approach},
  journal = {Communications in Mathematical Sciences},
  volume = {4},
  number = {3},
  pages = {575--590},
  year = {2006},
  doi = {10.4310/CMS.2006.v4.n3.a5}
}

@article{ZhangWangRong2011,
  author = {Zhang, Jing and Wang, Li-Lian and Rong, Zhijian},
  title = {A prolate-element method for nonlinear {PDE}s on the sphere},
  journal = {Journal of Scientific Computing},
  volume = {47},
  number = {1},
  pages = {73--92},
  year = {2011},
  doi = {10.1007/s10915-010-9421-y}
}

@article{DziukElliott2013,
  author  = {Dziuk, Gerhard and Elliott, Charles M.},
  title   = {Finite Element Methods for Surface {PDE}s},
  journal = {Acta Numerica},
  volume  = {22},
  pages   = {289--396},
  year    = {2013},
  doi     = {10.1017/S0962492913000056}
}

@article{DuJuTian2011,
  author  = {Du, Qiang and Ju, Lili and Tian, Li},
  title   = {Finite Element Approximation of the {Cahn--Hilliard} Equation on Surfaces},
  journal = {Computer Methods in Applied Mechanics and Engineering},
  volume  = {200},
  number  = {29--32},
  pages   = {2458--2470},
  year    = {2011},
  doi     = {10.1016/j.cma.2011.04.018}
}

@article{HakulaRasila2025,
  author  = {Hakula, Harri and Rasila, Antti},
  title   = {{Laplace--Beltrami} Equations and Numerical Conformal Mappings on Surfaces},
  journal = {SIAM Journal on Scientific Computing},
  volume  = {47},
  number  = {1},
  pages   = {A325--A342},
  year    = {2025},
  doi     = {10.1137/24M1656840}
}

@article{BaumgartHessWebb2003,
  author  = {Baumgart, Tobias and Hess, Samuel T. and Webb, Watt W.},
  title   = {Imaging Coexisting Fluid Domains in Biomembrane Models Coupling Curvature and Line Tension},
  journal = {Nature},
  volume  = {425},
  number  = {6960},
  pages   = {821--824},
  year    = {2003},
  doi     = {10.1038/nature02013}
}

@article{JulicherLipowsky1993,
  author  = {J{\"u}licher, Frank and Lipowsky, Reinhard},
  title   = {Domain-Induced Budding of Vesicles},
  journal = {Physical Review Letters},
  volume  = {70},
  number  = {19},
  pages   = {2964--2967},
  year    = {1993},
  doi     = {10.1103/PhysRevLett.70.2964}
}

@article{CallanJonesSorreBassereau2011,
  author  = {Callan-Jones, Andrew and Sorre, Benjamin and Bassereau, Patricia},
  title   = {Curvature-Driven Lipid Sorting in Biomembranes},
  journal = {Cold Spring Harbor Perspectives in Biology},
  volume  = {3},
  number  = {2},
  pages   = {a004648},
  year    = {2011},
  doi     = {10.1101/cshperspect.a004648}
}

@article{ParthasarathyYuGroves2006,
  author  = {Parthasarathy, Raghuveer and Yu, Cheng-han and Groves, Jay T.},
  title   = {Curvature-Modulated Phase Separation in Lipid Bilayer Membranes},
  journal = {Langmuir},
  volume  = {22},
  number  = {11},
  pages   = {5095--5099},
  year    = {2006},
  doi     = {10.1021/la060390o}
}

@article{MichaletBensimon1995,
  author  = {Michalet, Xavier and Bensimon, David},
  title   = {Vesicles of Toroidal Topology: Observed Morphology and Shape Transformations},
  journal = {Journal de Physique II France},
  volume  = {5},
  number  = {2},
  pages   = {263--287},
  year    = {1995},
  doi     = {10.1051/jp2:1995128}
}

@article{BouzarMenasMueller2015,
  author  = {Bouzar, Lila and M{\'e}nas, Ferhat and M{\"u}ller, Martin Michael},
  title   = {Toroidal Membrane Vesicles in Spherical Confinement},
  journal = {Physical Review E},
  volume  = {92},
  number  = {3},
  pages   = {032721},
  year    = {2015},
  doi     = {10.1103/PhysRevE.92.032721}
}

@article{Ratz2016Benchmark,
  author  = {R{\"a}tz, Andreas},
  title   = {A Benchmark for the Surface {Cahn--Hilliard} Equation},
  journal = {Applied Mathematics Letters},
  volume  = {56},
  pages   = {65--71},
  year    = {2016},
  doi     = {10.1016/j.aml.2015.12.008}
}

@article{XuTang2006,
  author  = {Xu, Chuanju and Tang, Tao},
  title   = {Stability Analysis of Large Time-Stepping Methods for Epitaxial Growth Models},
  journal = {SIAM Journal on Numerical Analysis},
  volume  = {44},
  number  = {4},
  pages   = {1759--1779},
  year    = {2006},
  doi     = {10.1137/050628143}
}

@article{ShenYang2010,
  author  = {Shen, Jie and Yang, Xiaofeng},
  title   = {Numerical approximations of {Allen--Cahn} and {Cahn--Hilliard} equations},
  journal = {Discrete and Continuous Dynamical Systems. Series A},
  volume  = {28},
  number  = {4},
  pages   = {1669--1691},
  year    = {2010},
  doi     = {10.3934/dcds.2010.28.1669}
}

@article{ShenXuYang2018,
  author  = {Shen, Jie and Xu, Jie and Yang, Jiang},
  title   = {The scalar auxiliary variable ({SAV}) approach for gradient flows},
  journal = {Journal of Computational Physics},
  volume  = {353},
  pages   = {407--416},
  year    = {2018},
  doi     = {10.1016/j.jcp.2017.10.021}
}

@article{JiangZhangZhao2022,
  author  = {Jiang, Maosheng and Zhang, Zengyan and Zhao, Jia},
  title   = {Improving the accuracy and consistency of the scalar auxiliary variable ({SAV}) method with relaxation},
  journal = {Journal of Computational Physics},
  volume  = {456},
  pages   = {110954},
  year    = {2022},
  doi     = {10.1016/j.jcp.2022.110954}
}

\end{document}